\documentclass[a4paper,11pt]{article}

\usepackage{amssymb}
\usepackage{amsmath}
\usepackage{amsthm}
\usepackage{amscd}
\usepackage[english]{babel}
\usepackage[all]{xypic}
\usepackage[all]{xy}
\usepackage{url}
\usepackage{graphicx}
\usepackage{hyperref}

\usepackage{comment}

\newtheorem{theo}{Theorem}[section]
\newtheorem*{theo*}{Theorem}
\newtheorem{prop}[theo]{Proposition}
\newtheorem{lem}[theo]{Lemma}

\newtheorem{defi}[theo]{Definition}
\newtheorem{rema}[theo]{Remark}

\def \kbar {{\bar k}}
\newcommand{\dem}[1] {\paragraph{Proof{#1}: }}

\def \Romannumeral #1 {\expandafter\uppercase\expandafter {\romannumeral #1} }
\def \br {{\rm{Br\,}}}

\def \T {{\mathcal T}}
\def \C {{\mathcal C}}

\def \nr {{\rm nr}}

\def \pic {{\rm {Pic\,}}}

\def \gal {{\rm{Gal}}}

\def \calo {{\cal O}}
\def \spec {{\rm{Spec\,}}}

\def \Z {{\bf Z}}
\def \Q {{\bf Q}}
\def \L {{\bf L}}
\def \F {{\bf F}}

\def \NN {{\bf N}}
\def \RR {{\bf R}}

\def \cok {{\rm{coker\,}}}

\def \G {{\bf G}_m}

\def\Sha{\cyrille X}

\def\A{{\bf A}}

\def\smallsquare{\vbox{\hrule\hbox{\vrule height 1 ex\kern 1 ex\vrule}\hrule}}
\def\enddem{\hfill\smallsquare\vskip 3mm}

\def \abstract{\paragraph{R\'esum\'e. }}
\def \abstractbis{\paragraph{Abstract. }}
\def \codams{\paragraph{AMS codes~: }}

\DeclareFontFamily{U}{wncy}{}
\DeclareFontShape{U}{wncy}{m}{n}{%
   <5>wncyr5%
   <6>wncyr6%
   <7>wncyr7%
   <8>wncyr8%
   <9>wncyr9%
   <10>wncyr10%
   <11>wncyr10%
   <12>wncyr6%
   <14>wncyr7%
   <17>wncyr8%
   <20>wncyr10%
   <25>wncyr10}{}
\DeclareMathAlphabet{\cyrille}{U}{wncy}{m}{n}
\def\Sha{\cyrille X}

\def \et {{\textup{\'et}}}

\def\Cone{{\rm Cone}}

\def\Gm{{\mathbf{G}_m}}

\author{Cyril Demarche and David Harari}

\title{Duality for complexes of tori over a global field of positive
 characteristic}

\begin{document}
\maketitle

\abstractbis{If $K$ is a number field, arithmetic duality theorems for tori and complexes of tori over $K$ are crucial to understand local-global principles for linear algebraic groups over $K$. When $K$ is a global field of positive characteristic, we prove similar arithmetic duality theorems, including a Poitou-Tate exact sequence for Galois hypercohomology of complexes of tori. One of the main ingredients is Artin-Mazur-Milne duality theorem for fppf cohomology of finite flat commutative group schemes.}

\codams{11E72, 11G20, 14F20, 14H25}

\section{Introduction}\label{intro}

Let $K$ be a global field of characteristic $p \geq 0$ and let $\A_K$ denotes the ring of ad\`eles of $K$. Let $G$ be a reductive group over $K$, and $X$ be a torsor under $G$. We are interested in rational points on $X$, and more precisely, on various local-global principles associated to $X$: does $X$ satisfy the Hasse principle, i.e. does $X(\A_K) \neq \emptyset$ imply $X(K) \neq \emptyset$? If not, can we explain the failure using the so-called Brauer-Manin obstruction to the Hasse principle? Assuming that $X(K) \neq \emptyset$, can we estimate the size of $X(K)$ by studying the so-called weak and strong approximation on $X$ (with a Brauer-Manin obstruction if necessary), i.e. the closure of the set $X(K)$ in the topological space $X(\A_K^S)$, where $S$ is a (not necessarily finite) set of places of $K$ and $\A_K^S$ is the ring of $S$-ad\`eles (with no components in $S$)?

The answer to those questions is known in the case where $K$ is a number field, see for instance \cite{San}, Corollaries 8.7 and 8.13, for the Hasse principle and weak approximation, and \cite{demlms}, Theorem 3.14 for strong approximation. Note that in the number field case, similar results are known for certain 
non principal 
homogeneous spaces of $G$ (see \cite{Bor} or \cite{BD}).

In the case of a global field of positive characteristic, the answer is known for semisimple simply connected groups (thanks to works by Harder, Kneser, Chernousov, Platonov, Prasad), but the general case is essentially open (see \cite{Ros}, Theorem 1.9 for some related results). One strategy to attack the remaining local-global questions is similar to one that worked for number fields: arithmetic duality theorems for tori, and abelianization of Galois cohomology (see for instance \cite{demphd} and \cite{demlms} for the case of strong approximation over number fields). Indeed, given a reductive group $G$ over a field $L$ (e.g. $L$ is a global or a local field), one can construct a complex of tori of length two $C := [T_1 \to T_2]$, together with "abelianization maps" $H^i(L,G) \to H^i(L,C)$ (cohomology sets here are Galois cohomology or hypercohomology sets), such that the cohomology sets of $G$ can be computed via the abelian cohomology groups of $C$ and the Galois cohomology of a semisimple simply connected group associated to $G$. The latter is well-understood when $L$ is a local or global field. 

Motivated by the discussion above, this paper deals with arithmetic duality theorems for complexes of tori over global fields $K$
of positive characteristic; 
in characteristic $0$, we also get refinements of previously known results.

The aforementioned applications to the arithmetic of reductive groups and homogeneous spaces will be given in a future paper.

The main object is a two-term complex $C := [T_1 \to T_2]$ of $K$-tori $T_1$ and $T_2$, and we are particularly interested in its Galois hypercohomology groups $H^i(K, C)$. The main result of the paper can be summarized as follows: we get Poitou-Tate exact sequences relating global Galois cohomology groups $H^i(K,C)$ and local ones $H^i(K_v,C)$ - for any place $v$ of $K$ - via the cohomology of the dual object $\widehat{C}$ of $C$. To be more precise, let us introduce some notation: $K$ is the function field of a smooth, projective, and geometrically
integral curve $X$ over a finite field $k$. Let $X^{(1)}$ denote the set of closed points in $X$. If $A$ is a discrete
abelian group, 
then $A^*$ is the Pontryagin dual of homomorphisms from $A$ to $\Q / \Z$, and $A_\wedge$ denotes the completion $A_{\wedge}:=\varprojlim_{n \in \NN^*} A/n$.

We can now state one of the main results in the paper (see Theorem \ref{spt1bis}):
\begin{theo*} 
Let $C := [T_1 \to T_2]$ be a two-term complex of $K$-tori, and $\widehat{C} := [\widehat{T_2} \to \widehat{T_1}]$ be the dual complex, where $\widehat{T}$ is the module of characters of a torus $T$ ($T_1$ and $\widehat{T}_2$ are in degree $-1$). Then there is an exact sequence 
\begin{equation} 
\begin{CD}
0 @>>> H^{-1}(K,C)_{\wedge} @>>> \left[\prod^{'}_{v \in X^{(1)}}
H^{-1}(K_v,C)\right]_{\wedge} @>>>
H^2(K,\widehat C)^* \cr
&&&&&& @VVV \cr
&& H^1(K,\widehat C)^* @<<< \left[\prod^{'} _{v \in X^{(1)}} H^{0}(K_v,C)\right]_{\wedge}
@<<< H^0(K,C)_{\wedge} \cr
&& @VVV \cr
&& H^1(K,C) @>>> \bigoplus_{v \in X^{(1)}} H^1(K_v,C) @>>> H^0(K,\widehat C)^* \cr
&&&&&& @VVV \cr
 0 @<<< H^{-1}(K,\widehat C)^* @<<< \bigoplus_{v \in X^{(1)}} H^2(K_v,C)
@<<< H^2(K,C) \cr
\end{CD}
\end{equation}
\end{theo*}

In the case where $K$ is a number field, we also recover a generalization of \cite{demphd}, Theorems 6.1 and 6.3. In the function field case, some partial results related to this exact sequence for one single torus can be deduced from section 6 in \cite{gonz1}.

The main ingredient to prove the Theorem above is the so-called Artin-Mazur-Milne duality theorem for the fppf cohomology of finite flat commutative group schemes over open subsets of $X$ (see \cite{MilADT}, Theorem III.8.2 and \cite{DemH}, Theorem 1.1).

The structure of the paper is the following: section \ref{one} extends the construction and properties given in \cite{DemH} of fppf cohomology with compact support to the case of bounded complexes of finite flat group schemes. General properties of \'etale cohomology of complexes of tori and of their dual complexes are given in section \ref{onebis}. Section \ref{two} deals with applications of Artin-Mazur-Milne duality theorem to various duality statements for the \'etale cohomology of complexes of tori over open subsets $U$ of $X$. In section \ref{three}, one deduces several Poitou-Tate exact sequences for Galois cohomology from the results of section \ref{two}.

\section{Compact support hypercohomology} \label{one}

Let $K$ be the function field of a smooth, projective, and geometrically
integral curve $X$ over a finite field $k$. Let $U$ be a non empty
Zariski open subset of $X$. Denote by $U^{(1)}$ the set of closed points 
of $U$.

\smallskip

Let $\mathcal{C}=(C_p)_{p \in \Z}$
be a bounded complex of fppf sheaves over $U$. In this text, we define the dual of $\mathcal{C}$ to be the Hom-complex $\widehat{\mathcal{C}}$ defined by 
\[ \widehat{\mathcal{C}} := \underline{\textup{Hom}}^\bullet(\mathcal{C}, \Gm[1]) \, ,\]
following the sign conventions in \cite{SP}, Tag 0A8H or in \cite{bkiac}, X.5.1.
Note that there is a functorial morphism of complexes
\[\textup{Tot}(\mathcal{C} \otimes \widehat{\mathcal{C}}) \to \Gm[1]\]
mapping an element $c \otimes \varphi \in C_p \otimes \textup{Hom}(C_q, A_n)$ to $0$ if $p \neq q$, and to $(-1)^{p (n-1)} \varphi(c) \in A_n$ if $p=q$. 

With those conventions, if $\mathcal{C}$ is concentrated in degree $0$, i.e. $\mathcal{C} = \mathcal{F}$ with $\mathcal{F}$ a fppf sheaf, then $\widehat{\mathcal{C}}$ is the same as the Cartier dual $\mathcal{F}^D := \underline{\textup{Hom}}(\mathcal{F}, \Gm)$ attached in degree $-1$, i.e. $\widehat{\mathcal{C}} = \mathcal{F}^D[1]$ and the above pairing coincides with the obvious pairing $\mathcal{F} \otimes \mathcal{F}^D[1] \to \Gm[1]$ with no extra-sign.

Note also that for any bounded complex $\mathcal{C}$, we have a natural isomorphism of complexes $\widehat{\mathcal{C}[1]} \xrightarrow{\sim} \widehat{\mathcal{C}}[-1]$, given by a sign $(-1)^{n+1}$ in degree $n$. And given a morphism $f : A \to B$ of bounded complexes, we have a natural isomorphism of complexes $\widehat{\textup{Cone}(f)} \xrightarrow{\sim} \textup{Cone}(\widehat{f})[-1]$ such that the following diagram commutes
\begin{equation} \label{cone and dual}
\xymatrix{
\widehat{A[1]} \ar[r] \ar[d]^\sim & \widehat{\textup{Cone}(f)} \ar[r] \ar[d]^\sim & \widehat{B} \ar[r]^{\widehat{f}} \ar[d]^= & \widehat{A} \ar[d]^= \\
\widehat{A}[-1] \ar[r] & \textup{Cone}(\widehat{f})[-1] \ar[r] & \widehat{B} \ar[r]^{\widehat{f}} & \widehat{A} \, .
}
\end{equation}

If $\mathcal N$ is a commutative
group scheme (over $U$ or over $K$), its Cartier dual is
denoted ${\mathcal N}^D$. The Pontryagin dual of a topological abelian group
$A$ (consisting of continuous homomorphism from $A$ to $\Q/\Z$) is
denoted $A^*$. Unless explicitely specified, the topology used for
sheaves (resp. complex of sheaves) and cohomology (resp.
hypercohomology) is the fppf topology. 

\smallskip

For each closed point $v$ of $X$, the completion of $K$ at $v$
is denoted by $K_v$: it is a local field of characteristic $p$ with
finite residue field $\F_v$ (observe the slight difference of notation with 
\cite{DemH}, where $K_v$ stands for the henselization and $\widehat{K_v}$
for the completion). Denote by $\calo_v$ the ring of integers of $K_v$.
For every fppf sheaf ${\mathcal F}$ 
over $U$ with generic fibre $F$, recall (\cite{DemH}, Prop~2.1)
the long exact sequence (where the piece of notation $v \not \in U$ means 
that we consider all closed points of $X-U$).

\begin{equation} \label{suppcomp}
...\to H^i_c(U,{\mathcal F}) \to H^i(U,{\mathcal F}) \to \bigoplus_{v \not 
\in U} H^i(K_v,F) \to H^{i+1}_c(U,\mathcal F) \to...
\end{equation}		
There is also a long exact sequence 
\begin{equation} \label{longex1}
...\to H^i_c(U,{\mathcal F}') \to H^i_c(U,{\mathcal F}) \to
H^i_c(U,{\mathcal F}'') \to H^{i+1}_c(U,\mathcal F ') \to...
\end{equation}
associated to every short exact sequence 
$$0 \to {\mathcal F}' \to {\mathcal F} \to {\mathcal F}'' \to 0$$
of fppf sheaves.

\smallskip

Let us now extend the construction of the groups $H^i_c(U,...)$ and \cite{DemH}, Prop~2.1 to the case of bounded complexes.
Let $\mathcal{C} := [\dots \to \mathcal{F}_i \to \mathcal{F}_{i+1} \to \dots]$ be a bounded complex of fppf sheaves over $U$.
Let $\mathcal{C} \to I^\bullet(\mathcal{C})$ be an injective resolution of the complex $\mathcal{C}$, in the sense of \cite{SP}, Tag 013K.
Following \cite{DemH}, section 2, let $Z := X \setminus U$ and
$Z' := \coprod_{v \in Z} \spec({K_v}) \xrightarrow{i} U$.
 Denote by $\mathcal{C}_v$ and $I^\bullet(\mathcal{C})_v$ their respective 
pullbacks to $\spec K_v$, for $v \notin U$.

We now define $\Gamma_c(U, I^\bullet(\mathcal{C}))$ to be the following object in the category of complexes of abelian groups:
$$\Gamma_c(U, I^\bullet(\mathcal{C})) := \Cone\left(\Gamma(U, I^\bullet(\mathcal{C})) \to \Gamma(Z', i^* I^\bullet(\mathcal{C}))  \right)[-1] \, ,$$
and $H^r_c(U, \mathcal{C}) := H^r(\Gamma_c(U, I^\bullet(\mathcal{C})))$. We will also denote by $R\Gamma_c(U,\mathcal{C})$ the complex $\Gamma_c(U, I^\bullet(\mathcal{C}))$. 
Similarly, one can define, for any closed point $v \in X$, complexes $\Gamma_v(\mathcal{O}_v, \mathcal{C})$ computing fppf cohomology groups $H^r_v(\mathcal{O}_v, \mathcal{C})$ over $\spec \mathcal{O}_v$ with support in the closed point, as in \cite{DemH}, before Lemma 2.6. 

As in \cite{DemH}, similar definitions could be made when $K$ is a number field (taking into account the real places), but in this article we will focus on the function field case. However, we will make remarks regularly throughout the text explaining similarities and differences appearing in the number field case.

We will need the analogue of \cite{DemH}, Prop~2.1 and Prop~2.12 for bounded complexes $\mathcal{C}$: by construction, the first two points of 
loc.~cit. Prop.~2.1 (i.e. exact sequence \eqref{suppcomp} and
(\ref{longex1})) still hold for bounded complexes.
\begin{prop} \label{prop 2.1.3 and 2.12}
Let $\mathcal{C}$ be a bounded complex of flat affine commutative group schemes of finite type over $U$, and let $V \subset U$ be a non empty open subset. 
\begin{enumerate}
	\item There is a canonical commutative diagram of abelian groups:
\begin{changemargin}{-4cm}{4cm}
\[
\xymatrix{
& & \bigoplus_{v \notin V} H^{r-1}({K_v},\mathcal{C}) \ar[d] & \bigoplus_{v \notin U} H^{r-1}({K_v}, \mathcal{C}) \ar[l]_-{i_2} \ar[d] & &\\
\dots \ar[r] & \bigoplus_{v \in U \setminus V} H^{r-1}({\mathcal{O}_v}, \mathcal{C})
 \ar[r] \ar[ru]^-{i_1} & H^r_c(V, \mathcal{C}) \ar[r] \ar[d] & H^r_c(U, \mathcal{C}) \ar[r] \ar[d] & \bigoplus_{v \in U \setminus V} H^{r}({\mathcal{O}_v}, \mathcal{C}) \ar[r] & \dots \\
& & H^r(V,\mathcal{C}) \ar[d]& H^r(U,\mathcal{C}) \ar[d] \ar[l]_-{\textup{Res}} \ar[ru] & & \\
& & \bigoplus_{v \notin V} H^r({K_v},\mathcal{C}) \ar[r]^-{\pi} & \bigoplus_{v \notin U} H^r({K_v}, \mathcal{C}) \, , & &
}
\]
\end{changemargin}
where the long row and the columns are exact.
	\item Let $V \subset U$ be a non empty open subset. Then there is 
an exact sequence 
\[...\to \bigoplus_{v \in U \setminus V} H^r_v(\calo_v,\mathcal C) \to 
H^r(U,\mathcal C) \to H^r(V,\mathcal C) \to \bigoplus_{v \in U \setminus V}
H^{r+1}_v(\calo_v,\mathcal C) \to...\]
\end{enumerate}
\end{prop}

\dem{}
We follow the proofs of \cite{DemH}, Prop~2.1 3 and of Prop~2.12. 

Easy d\'evissages imply that \cite{DemH}, Lemmas~2.6 holds with $\mathcal{F}$ replaced by a bounded complex of flat commutative group schemes of finite type and Lemma~2.9 holds for bounded complexes of fppf sheaves. Likewise 
Lemma~2.10 holds for bounded complexes of \'etale sheaves
or of smooth commutative group schemes.
Therefore, one can copy the proofs of \cite{DemH},
Prop~2.1 (3) and Prop~2.12 to get the required Proposition.
\enddem

\begin{lem} \label{uvlem}
Let ${\mathcal C}$ be a bounded complex of flat commutative group schemes 
of finite type over $U$ with generic fibre $C$ over $K$. Let $i$ be an integer. 
For each $v \in U^{(1)}$, denote by $H^i_{\nr}(K_v,C)$ the image of
$H^i(\calo_v,\mathcal C)$ in $H^i(K_v,C)$. 
Let $V \subset U$ be a non empty Zariski open subset.
Then there is an exact sequence 
$$H^i(U,\mathcal C) \to \prod_{v \not \in U} H^i(K_v,C) \times
\prod_{v \in U-V} H^i_{\nr}(K_v,C) \to H^{i+1}_c(V,\mathcal C).$$
\end{lem}

\dem{} There is a commutative diagram such that the second line and the 
left column are exact (by (\ref{suppcomp}) and Prop.~\ref{prop 2.1.3 and 2.12} 2.):

$$
\begin{CD}
H^i(U,\mathcal C) @>>> \prod_{v \not \in U} H^i(K_v,C) \times
\prod_{v \in U-V} H^i_{\nr}(K_v,C) \cr
@VVV @VVjV \cr 
H^i(V,\mathcal C) @>>> \prod_{v \not \in U} H^i(K_v,C) \times 
\prod_{v \in U-V} H^i(K_v,C) @>>> H^{i+1}_c(V,\mathcal C) \cr
@VVV @VVV \cr 
\prod_{v \in U-V} H^{i+1}_v(\calo_v,\mathcal C) @>>=> \prod_{v \in U-V} H^{i+1}_v(\calo_v,\mathcal C)
\end{CD}
$$

For $v \in U-V$, the localization exact sequence (cf. proof of 
Proposition \ref{prop 2.1.3 and 2.12}, 2.)
$$H^i(\calo_v,\mathcal C) \to H^i(K_v,C) \to H^{i+1}_v(\calo_v,\mathcal C)$$
yields that 
the second column is a complex. Since $j$ is injective by definition, the 
required exact sequence follows by diagram chasing.

\enddem

For a complex $\mathcal{C}$ of finite flat group schemes, 
let us now endow the groups $H^*_c(U, \mathcal{C})$ with a natural topology, compatible with the one defined in \cite{DemH} in the case where $\mathcal{C}$ is a finite flat group scheme.

Let $F$ be a local field (that is: a field complete for a discrete valuation 
with finite residue field) and let $\mathcal{C} := [C_r \xrightarrow{f_r} C_{r+1} \to \dots \to C_s]$ be a bounded complex of finite commutative group schemes over $\spec F$, with $C_i$ in degree $i$. We assume that $F$ is of positive 
characteristic $p$ (if $F$ is $p$-adic, then all 
groups $H^r(F,\mathcal C)$ are finite by \cite{MilADT}, Cor. I.2.3). 

\begin{defi}
{\rm A morphism $f : G_1 \to G_2$
of topological groups is {\it strict} if it is continuous, and
the restriction $f : G_1 \to f(G_1)$ is
an open map (where the topology on $f(G_1)$ is induced by $G_2$). This is
equivalent to saying that $f$ induces an isomorphism of the topological
quotient $G_1/\ker f$ with the topological subspace $f(G_1) \subset G_2$.
}
\end{defi}

Let $A := \ker(f_r)$, and let $\overline{\mathcal{C}} := \textup{Cone}(A[-r] \xrightarrow{j} \mathcal{C})$. Then there is an exact triangle
\begin{equation} \label{devissage complexe}
A[-r] \xrightarrow{j} \mathcal{C} \xrightarrow{i} \overline{\mathcal{C}} \xrightarrow{p} A[1-r] \, .
\end{equation}
In addition, we have a natural quasi-isomorphism $\varphi : \overline{\mathcal{C}} \to \mathcal{C}'$, where $\mathcal{C}' := [C_{r+1}/{\rm Im} \, (f_r) \to C_{r+2} \to \dots \to C_s]$ has a smaller length than $\mathcal{C}$.

There is an alternative d\'evissage for the complex $\mathcal{C}$, given by
the exact triangle:
\begin{equation} \label{devissage complexe bis}
\widetilde{\mathcal{C}} \xrightarrow{i'} \mathcal{C} \xrightarrow{p} C_r[-r] \xrightarrow{\partial} \widetilde{\mathcal{C}}[1] \, ,
\end{equation}
where $\widetilde{\mathcal{C}} := [C_{r+1} \xrightarrow{f_{r+1}} C_{r+2} \to \dots \to C_s]$. 

Recall that for a finite and commutative $F$-group scheme $N$, the fppf groups 
$H^i(F,N)$ are finite if $i \neq 1$ (\cite{MilADT}, \S III.6) and they are
equipped with a locally compact topology for $i=1$ by \cite{Ces}.
By induction on the length of $\mathcal{C}$, one deduces that if $C_i = 0$, then $H^{i+1}(F, \mathcal{C})$ is finite. In particular, with the previous notation, we get that $H^i(F, \mathcal{C})$ is finite if $i \leq r$ or $i \geq s+2$.

We now define a natural topology on $H^i(F, \mathcal{C})$ by induction on the length of $\mathcal{C}$, so that any morphism of such complexes induces a strict map between hypercohomology groups. Using the d\'evissages given by \eqref{devissage complexe} and \eqref{devissage complexe bis}, one gets the following exact sequences 
\begin{equation} \label{devissage cohomo}
H^{i-1}(F, \mathcal{C}') \to H^{i-r}(F,A) \xrightarrow{f} H^i(F, \mathcal{C}) \xrightarrow{g} H^i(F, \mathcal{C}') \to H^{i+1-r}(F,A) \to \dots
\end{equation}
and 
\begin{equation} \label{devissage cohomo bis}
H^{i-r-1}(F, C_r) \to H^{i}(F,\widetilde{\mathcal{C}}) \xrightarrow{f'} H^i(F, \mathcal{C}) \xrightarrow{g'} H^{i-r}(F, C_r) \to H^{i+1}(F,\widetilde{\mathcal{C}}) \to \dots \, ,
\end{equation}
where all the groups except $H^i(F, \mathcal{C})$ are endowed with a
natural topology via the induction hypothesis (observe that $C'_r=0$). 
\begin{itemize}
    \item Assume that $i = r+1$. Equip ${\rm Im } \, f\cong H^{i-r}(F,A)/
{\rm Im} \, H^{i-1}(F, \mathcal{C}')$ with the quotient topology. 
Then the two rightmost groups in exact sequence \eqref{devissage cohomo} are finite, and we can endow $H^{r+1}(F, \mathcal{C})$ with the topology such that 
${\rm Im} \, f$ is an open subgroup (see \cite{DemH}, beginning of section 3). 
In other words it is the 
finest topology such that $f$ is continuous. Then $f$ and $g$ are strict.
    \item Assume that $i \neq r+1$. Then $H^{i-r}(F,C_r)$ is finite (and discrete), and using exact sequence \eqref{devissage cohomo bis} one can similarly 
endow $H^i(F, \mathcal{C})$ with the finest topology making $f'$ continuous. Then both maps $f'$ and $g'$ are strict.
\end{itemize}
By construction, this topology is functorial in $\mathcal{C}$, i.e. if $\mathcal{C}_1 \to \mathcal{C}_2$ is a morphism of complexes, then the induced map $H^i(F, \mathcal{C}_1) \to H^i(F, \mathcal{C}_2)$ is strict. In addition, given a quasi-isomorphism $\mathcal{C}_1 \to \mathcal{C}_2$, the induced morphism on cohomology groups is a homeomorphism.

Let us now deal with the topology on the groups $H^*_c(U, \mathcal{C})$, where $\mathcal{C}$ is a complex of finite flat commutative group schemes defined over $U$. Recall that we have an exact sequence analogous to \eqref{suppcomp}:
\begin{equation} \label{long exact seq compact supp}
H^{i-1}(U,\mathcal{C}) \to \bigoplus_{v \notin U} H^{i-1}(K_v,\mathcal{C}) \to H^{i}_c(U,\mathcal{C}) \to H^{i}(U,\mathcal{C}) \, .
\end{equation}
We endow the groups $H^{i}(U,\mathcal{C})$ with the discrete topology, and the groups $H^{i-1}(K_v,\mathcal{C})$ with the topology defined above.

\begin{lem} \label{closedimagecomplex}
For all $i$, the map $H^i(U, \mathcal{C}) \to \bigoplus_{v \notin U} H^{i}(K_v,\mathcal{C})$ has discrete image.
\end{lem}

\dem{}
We prove the result by induction on the length of $\mathcal{C}$, using the d\'evissages induced by the exact triangles \eqref{devissage complexe} and \eqref{devissage complexe bis}.
\begin{itemize}
    \item Assume that $i = r+1$. Using the exact triangle \eqref{devissage complexe bis}, we get the following commutative diagram of long exact sequences of topological groups (where all maps are strict):
    \[
    \xymatrix{
    H^{r+1}(U, \widetilde{\mathcal{C}}) \ar[r] \ar[d] & H^{r+1}(U, \mathcal{C}) \ar[r] \ar[d] & H^1(U, C_r) \ar[d] \\
    \bigoplus_{v \notin U} H^{r+1}(K_v, \widetilde{\mathcal{C}}) \ar[r] & \bigoplus_{v \notin U} H^{r+1}(K_v, \mathcal{C}) \ar[r] & \bigoplus_{v \notin U} H^1(K_v, C_r) \, .
    }
    \]
    The groups on the left hand side are finite, and by \cite{cesnalms} Lemma 2.7, the image of the right hand side map is discrete in $\bigoplus_{v \notin U} H^1(K_v, C_r)$. Since $\bigoplus_{v \notin U} H^{r+1}(K_v, \mathcal{C})$ is Hausdorff, an easy topological argument implies that the image of the central vertical map is discrete.
    \item Assume that $i \neq r+1$. Using the exact triangle \eqref{devissage complexe}, we get the following commutative diagram of long exact sequences of topological groups:
    \[
    \xymatrix{
    H^{i-r}(U, A) \ar[r] \ar[d] & H^{i}(U, \mathcal{C}) \ar[r] \ar[d] & H^i(U, \mathcal{C}') \ar[d] \\
    \bigoplus_{v \notin U} H^{i-r}(K_v, A) \ar[r] & \bigoplus_{v \notin U} H^{i}(K_v, \mathcal{C}) \ar[r] & \bigoplus_{v \notin U} H^i(K_v, \mathcal{C}') \, .
    }
    \]
    The groups on the left hand side are finite, and by induction on the length of the complex, the image of the right hand side map is discrete in $\bigoplus_{v \notin U} H^1(K_v, \mathcal{C}')$. A similar topological argument as before implies that the central vertical map is discrete.
\end{itemize}
\enddem

As a consequence of this Lemma, one can endow $H^i_c(U, \mathcal{C})$ with the following topology: we put the quotient topology on the group $\bigoplus_{v \notin U} H^{i}(K_v,\mathcal{C}) / {\rm Im} \, H^i(U, \mathcal{C})$ (this topology is Hausdorff),
 and since $H^i(U, \mathcal{C})$ is discrete, there is a unique topology on $H^i_c(U, \mathcal{C})$ so that the maps in the exact sequence \eqref{long exact seq compact supp} are strict.

\begin{lem} \label{lem topo profinite}
The topological group $H^i_c(U, \mathcal{C})$ is profinite.
\end{lem}

\dem{}
We prove this Lemma by induction on the length of the complex $\mathcal{C}$.
By \cite{DemH}, Proposition 3.5, the Lemma is proven when $\mathcal{C}$ is a complex of length one, i.e. concentrated in one given degree.

Given a complex $\mathcal{C}$, consider the previous d\'evissages:
\begin{itemize}
    \item Assume that $i = r+1$. Then exact sequence \eqref{devissage cohomo} implies that the group $H^{r+1}_c(U, \mathcal{C})$ is an extension (the maps 
being strict) of a (discrete) finite group by a profinite group (which is a quotient of $H^1_c(U, A)$ by a closed subgroup), hence $H^{r+1}_c(U, \mathcal{C})$ is profinite.
    \item Assume that $i \neq r+1$. Then exact sequence \eqref{devissage cohomo bis} implies that $H^i_c(U, \mathcal{C})$ is an extension of a finite (discrete) group by a profinite group (which is a quotient of $H^i_c(U, \widetilde{\mathcal{C}})$ by a closed subgroup), hence it is profinite.
\end{itemize}
\enddem

\section{Cohomology of tori and short complexes of tori} \label{onebis}

Let $U$ be a non empty
Zariski open subset of $X$. Recall that for every {\it $U$-torus}
$\T$ (in the sense of \cite{sga3}, IX, D\'ef.~1.3),
there is a finite \'etale covering
(that can be taken to be connected and
Galois) $V$ of $U$ such that $\T_V:=\T \times_U V$
is {\it split}, that is: isomorphic to some power $\G ^r$ ($r \in \NN$) of the
multiplicative group  (\cite{sga3}, X, Th.~5.16). The {\it group of characters}
$\widehat \T$ of $\T$ is a $U$-group scheme locally isomorphic to
$\Z^r$ for the \'etale topology, namely it is a torsion-free and finite
type $\gal(V/U)$-module.

\smallskip

Given a complex of $U$-tori $\C=[\T_1 \stackrel{\rho}{\to} \T_2]$ 
with generic fibre $C=[T_1 \stackrel{\rho}{\to} T_2]$ over $K$,
where by convention $\T_1$ is in degree $-1$ and $\T_2$ in degree $0$, we
can apply the construction of section \ref{one}.
Namely we have dual complexes $\widehat {\mathcal C}
=[\widehat {\mathcal T}_2 \to \widehat {\mathcal T}_1]$ and
$\widehat C=[\widehat T_2 \stackrel{\hat \rho}{\to}
\widehat T_1]$ (concentrated in degrees $-1$ and $0$), which are
respectively defined over $U$ and over $K$.
Denote by $S$ the finite set $X-U$ and by $G_S=\pi_1^{\et}(U)$ the \'etale fundamental
group of $U$, which is the
Galois group of the maximal separable
field extension $K_S$ of $K$ unramified outside $S$; then
each $\widehat \T_i$ ($i=1,2$) can be viewed as a discrete $G_S$-module.

\smallskip

Recall that fppf and \'etale
cohomology coincide for sheaves represented by smooth group schemes
(\cite{MilEC}, section III.3) like a torus $\T$, its group
of characters $\widehat \T$, or finite flat group schemes of order
prime to $p$. In particular (by \cite{MilEC}, Lemma III.1.16) we have
for every integer $i$:
$$\varinjlim_U H^i(U,\C) \cong H^i(K,C)$$
(where the limit is over all non empty Zariski open subsets $U$ of $X$),
and likewise for the complex $\widehat \C$.

\smallskip

For such 2-term complexes, the pairings of section \ref{one} 
can be made explicit (see \cite{demphd}, section 2; note that the sign conventions are slightly different here), and give maps
$${\mathcal C} \otimes^{{\bf L}}
\widehat {\mathcal C} \to \G[1]; \quad C \otimes^{\bf L} \widehat C \to \G[1]
$$
in the bounded derived category ${\mathcal D}^b(U)$ (resp. ${\mathcal D}^b(\spec K)$)
of fppf sheaves over $U$ (resp. over $\spec K$).
In the case $T_1=0$ or $T_2=0$, we recover (up to shift) the
classical pairings $\T \otimes \widehat \T \to \G$ and $T \otimes
\widehat T \to \G$ associated to one single torus $\T$.
We also have for each positive integer $n$ the $n$-adic realizations
$$T_{\Z/n}(\C):=H^0(\C[-1] \otimes^{\bf L} \Z/n); \quad
T_{\Z/n}(\widehat \C):=H^0(\widehat \C[-1] \otimes^{\bf L} \Z/n)$$
and likewise for $C$ and $\widehat C$.
The fppf sheaf $T_{\Z/n}(\C)$ is representable by a finite group
scheme of multiplicative type over $U$ (in the sense of \cite{sga3}, IX,
D\'ef.~1.1) with Cartier dual
$T_{\Z/n}(\widehat \C)$, and similarly for $T_{\Z/n}(C)$ and
$T_{\Z/n}(\widehat C)$ over $K$. Besides we have exact
triangles (\cite{demphd}, Lemme~2.3), where for every abelian group (or group
scheme) $A$, the piece of notation
$_ n A$ stands for the $n$-torsion subgroup of $A$:
\begin{equation} \label{exactriang1}
_n(\ker \rho)[2] \to \C \otimes^{\bf L} \Z/n \to
T_{\Z/n}(\C)[1] \to \, _n(\ker \rho)[3]
\end{equation}
and
\begin{equation} \label{exactriang2}
T_{\Z/n}(\widehat \C)[1] \to {\widehat \C} \otimes^{\bf L} \Z/n \to
 \widehat{_n (\ker \rho)} \to T_{\Z/n}(\widehat \C)[2]
\end{equation}
in ${\mathcal D}^b(U)$, and similar triangles for $C$,
$\widehat C$ in ${\mathcal D}^b(\spec K)$.

Note also that the objects $\mathcal{C} \otimes^{\bf L} \Z / n$ and $\widehat{\mathcal{C}} \otimes^{\L} \Z / n$ in the derived category have canonical representatives as complexes of fppf sheaves given by $\mathcal{C} \otimes^{\bf L} \Z / n = \textup{Tot}(\mathcal{C} \otimes [\Z \xrightarrow{n} \Z])$ and $\widehat{\mathcal{C}} \otimes^{\L} \Z / n = \textup{Tot}(\widehat{\mathcal{C}} \otimes \Z / n)$.

\smallskip

We also have an exact triangle in ${\mathcal D}^b(U)$:
\begin{equation} \label{devistm}
(\ker \rho)[1] \to \mathcal C \to \cok \rho \to (\ker \rho)[2],
\end{equation}
where $\cok \rho$ is a torus and $M:=\ker \rho$ is a group of multiplicative
type, and the dual exact triangle
\begin{equation} \label{devistmbis}
(\widehat{\cok \rho})[1] \to \widehat \C \to \widehat{\ker \rho}
\to (\widehat{\cok \rho})[2].
\end{equation}

For every integer $i$, there are exact sequences
\begin{equation} \label{devis1}
...\to H^i(U,{\mathcal T}_1) \to H^i(U,{\mathcal T}_2) \to
H^i(U,{\mathcal C}) \to H^{i+1}(U,{\mathcal T}_1) \to...
\end{equation}
\begin{equation} \label{devis2}
...\to H^i(U,\widehat {\mathcal T}_2) \to H^i(U,\widehat {\mathcal T}_1) \to
H^i(U,\widehat {\mathcal C}) \to H^{i+1}(U,\widehat {\mathcal T}_2) \to...
\end{equation}
and we also have similar exact sequences for the compact support fppf 
cohomology groups $H^i_c(U,...)$.

\begin{lem} \label{fini1}

\smallskip

a) Let ${\mathcal T}$ be a torus over $U$. Then 
$H^0(U,{\mathcal T})$ and $H^1(U,\mathcal T)$
are of finite type. If $U \neq X$, then 
$H^1(U,{\mathcal T})$ is finite.

\smallskip

b) Let ${\mathcal N}$ be a finite group scheme of multiplicative type over 
$U$. Then $H^1(U,{\mathcal N})$ and $H^2_c(U,{\mathcal N}^D)$ are finite. 
\end{lem}

\dem{} a) Let $V$ be an \'etale and finite connected 
Galois covering of $U$ such that ${\mathcal T}_V:={\mathcal T} \times_U V$
is isomorphic to $\G^r$ for some non negative integer $r$.  
The group ${\mathcal T}(V) \simeq
H^0(V,\G)^r$ is of finite type by Dirichlet's theorem on units. Therefore 
$H^0(U,{\mathcal T}) \subset H^0(V,{\mathcal T})={\mathcal T}(V)$
is also of finite type.

\smallskip

Set $G=\gal(V/U)$. Then Hochschild-Serre spectral sequence provides an 
exact sequence 
$$0 \to H^1(G,{\mathcal T}(V)) \to H^1(U,{\mathcal T})
 \to H^1(V,{\mathcal T}_V).$$
Since ${\mathcal T}_V \cong \G^r$, the group 
$H^1(V,{\mathcal T}_V) \cong (\pic V)^r$ is of finite type 
(resp. finite if $U \neq X$) by finiteness of the ideal class group
of a global field. 
As ${\mathcal T}(V)$ is of finite type and $G$ is finite, the group
$H^1(G,{\mathcal T}(V))$ is finite by \cite{col}, Chap. VIII, Cor 2.
Thus $H^1(U,{\mathcal T})$ is of 
finite type (resp. finite if $U \neq X$) as well. 

\smallskip

b) We can assume (by (\cite{MilADT}, Lemma III.8.9
and \cite{DemH}, Cor.~4.8) that $U \neq X$. 
Since every finitely generated Galois module is a quotient of a 
torsion-free and finitely generated Galois module, the 
assumption that ${\mathcal N}$ is of multiplicative type implies
(by \cite{sga3}, X, Prop~1.1) 
that there is an exact sequence of $U$-group schemes
$$0 \to {\mathcal N} \to {\mathcal T}_1 \to {\mathcal T}_2 \to 0,$$
where ${\mathcal T}_1$ and ${\mathcal T}_2$ are $U$-tori. Therefore there is 
an exact sequence of abelian groups 
$$ H^0(U,{\mathcal T}_2) \to H^1(U, {\mathcal N}) \to H^1(U, {\mathcal T}_1).$$
By a), we know that $H^1(U, {\mathcal T}_1)$ is finite. Let $n$ be the order
of ${\mathcal N}$; then the map $H^0(U,{\mathcal T}_2) \to H^1(U, {\mathcal N})$
factorizes through a map $H^0(U,{\mathcal T}_2)/n \to H^1(U, {\mathcal N})$. 
But $H^0(U,{\mathcal T}_2)/n$ is finite because $H^0(U,{\mathcal T}_2)$
is of finite type by a). Finally $H^1(U, {\mathcal N})$ is finite. The 
finiteness of $H^2_c(U,N^D)$ follows by Artin-Mazur-Milne
duality (\cite{DemH}, Th.~1.1).

\enddem

\begin{rema} \label{fini2}
{\rm By d\'evissage, the finiteness of $H^1(U, {\mathcal N})$ holds for 
a (not necessarily finite)
group of multiplicative type $\mathcal N$
because such a group is an extension of 
a finite group by a torus. Recall also (\cite{MilADT}, Lemma III.8.9 
and \cite{DemH}, Cor.~4.8) 
that for every finite and flat commutative 
group scheme ${\mathcal N}$ over $U$, the groups $H^i(U,{\mathcal N})$ and 
$H^{3-i}_c(U,{\mathcal N})$ are finite if $i \neq 1$ or if $U=X$, and 
also if $p$ does not divide the order of ${\mathcal N}$ (by \cite{DemH}, 
Prop~2.1., 4. and \cite{MilADT}, Th.~II.3.1).
Besides these groups are trivial if $i \geq 4$ (this is part of \cite{DemH}, 
Th.~1.1). 
}
\end{rema}

For a fppf sheaf (or a bounded complex of fppf sheaves)
${\mathcal F}$ on $U$ with generic fibre $F$, we set (cf. exact
sequence (\ref{suppcomp}))
$$D^i(U,{\mathcal F})={\rm Ker} \, [H^i(U, {\mathcal F}) \to 
\bigoplus_{v \not \in U} H^i(K_v,F)]={\rm Im} \, 
[ H^i_c(U,{\mathcal F}) \to H^i(U,{\mathcal F})].$$

\begin{lem} \label{bhn}
We have $D^2(U,\G)=0$. 
\end{lem}

\dem{} Let $\F_v$ be the residue field of $K_v$.
By \cite{MilADT}, Prop~II.1.1 b),
we have $$H^2(\calo_v,\G)=\br \calo_v \cong \br \F_v=0$$ 
because $\F_v$ is finite (\cite{nsw}, Prop. VI.3.5).
The Brauer group $\br U$ of $U$ injects into $\br K$
(\cite{MilEC}, Cor. IV.2.6). Now every element of
$D^2(U,\G) \subset \br U \subset \br K$ has trivial restriction to $\br K_v$ 
for all places $v$ of $K$, hence it is trivial by Brauer-Hasse-Noether 
Theorem (\cite{nsw}, Th.~VIII.1.17).

\enddem

\begin{lem} \label{h1lem}

Let $\mathcal T$ be a $U$-torus with generic fibre $T$. 

\smallskip 

a) The group $H^1(U,\widehat{\mathcal T})$ is finite;
the groups $H^0(U,\widehat{\mathcal T})$ and $H^0_c(U,\widehat{\mathcal T})$ 
are of finite type and torsion-free. The group $H^1_c(U,\widehat{\mathcal T})$
is of finite type.

\smallskip 

b) The group $H^2_c(U,\mathcal T)$ is finite. In particular, if $U=X$, then 
$H^2(X,\T)$ is finite. 

\smallskip

c) Assume $U \neq X$. Then $H^2_c(U, \widehat {\mathcal T})$ is finite.

\smallskip

d) Assume $i \geq 4$. Then $H^i(U,\widehat {\mathcal T})=H^i_c(U,\widehat
{\mathcal T})=0$. If $U \neq X$, 
then $H^3(U,\widehat {\mathcal T})=0$. 

\smallskip

e) If $U=X$, then $H^3(U,\widehat \T)=H^3(X,\widehat \T)$ is finite.

\end{lem}

\dem{} a) Let $V$ be a finite connected Galois \'etale covering of $U$ such that
${\mathcal T} \times_U V$ is split. Let $L$ be the function field 
of $V$ and set $G=\gal(L/K)$. We have
$H^1(V,\Z)=0$ because $H^1(V,\Z)$ injects into $H^1(L, \Z)$
by Leray spectral sequence.  
Therefore we have $H^1(V,\widehat{\mathcal T})=0$,
hence the group $H^1(U,\widehat {\mathcal T})$ identifies
(by Hochschild-Serre spectral sequence) to a subgroup of
$H^1(G,\widehat T)$, which is finite because $\widehat T$ is a
$G$-module of finite type. The assertion about $H^0(U,\widehat{\mathcal T})$
(which is a subgroup of $H^0(K,\widehat T)$) and $H^0_c(U,\widehat{\mathcal T}) 
\subset H^0(U,\widehat{\mathcal T})$ is obvious (we even have 
$H^0_c(U,\widehat{\mathcal T})=0$ if $U \neq X$). Also the exact sequence 
$$\bigoplus_{v \not \in U} H^0(K_v,\widehat T) \to H^1_c(U,\widehat \T) \to 
H^1(U,\widehat \T)$$ shows shat $H^1_c(U,\widehat \T)$ is of finite type.

\smallskip

b) By (\ref{suppcomp}), there is an exact sequence 
$$\bigoplus_{v \not \in U} H^1(K_v,T) \to H^2_c(U, \mathcal T)
 \to D^2(U,\mathcal T) \to 0 . $$
The groups $H^1(K_v,T)$ are finite (\cite{MilADT}, Cor. I.2.3). 
Since we have $D^2(V,\G)=0$ by Lemma~\ref{bhn},
a restriction-corestriction argument 
shows that $D^2(U,\mathcal T)$ is a subgroup of $_n H^2(U,\mathcal T)$, 
where $n=\# \gal(V/U)$. 
By Kummer sequence in the fppf topology
$$ 0 \to \, _n \mathcal T \to \mathcal T \stackrel{.n}{\to} \mathcal T \to 0,$$
the group $_n H^2(U,\mathcal T)$ is a quotient of 
$H^2(U,\, _n \mathcal T)$, which is finite (even if $p$ divides $n$, 
cf. Remark~\ref{fini2}). Thus $H^2_c(U,\mathcal T)$ is finite.

\smallskip

c) Let $n >0$. Using the exact sequence in the fppf topology
\begin{equation} \label{dualkummer}
0 \to \widehat{\mathcal T} \stackrel{.n}{\to} \widehat{\mathcal T} \to
\widehat{\mathcal T}/n \to 0 ,
\end{equation}
we see that $_n H^2_c(U,\widehat{\mathcal T})$ is a quotient of 
$H^1_c(U,\widehat{\mathcal T}/n)$, which is finite (see Remark~\ref{fini2}). 
It is therefore sufficient to show that $H^2_c(U,\widehat{\mathcal T})$ is 
of finite exponent, and by a restriction-corestriction argument, we 
reduce to the case $\widehat{\mathcal T}=\Z$. As $H^1(K_v,\Z)=0$, 
we have 
$$H^2_c(U,\Z)={\rm Ker} [H^2(U,\Z) \to \bigoplus_{v \not \in U} H^2(K_v,\Z)].$$
By \cite{MilADT}, Lemma II.2.10, this yields 
$$H^2_c(U,\Z) \cong {\rm Ker} [H^1(U,\Q/\Z) \to
\bigoplus_{v \not \in U} H^1(K_v,\Q/\Z)],$$ hence 
$$H^2_c(U,\Z) \cong 
{\rm Ker} [H^1(G_S,\Q/\Z) \to \bigoplus_{v \not \in U} H^1(K_v,\Q/\Z)].$$
Therefore $H^2_c(U,\Z)$ is a subgroup of $H^1(\gal(L/K),\Q/\Z)$, where 
$L \subset K_S$
is the maximal abelian extension of $K$ that is unramified outside $S$ 
and totally decomposed at every $v \in S=X-U$.  
The group $\gal(L/K)$ is isomorphic to $\pic U$ by class field 
theory, which implies that it is finite because $U \neq X$. Hence 
$H^2_c(U,\Z)$ is finite, which proves the lemma. 

\smallskip

d) For $i \geq 4$, the groups $H^i_c(U,\widehat \T)$ and $H^i(U, \widehat \T)$
coincide thanks to exact sequence (\ref{suppcomp}) because 
the local field $K_v$ is of strict cohomological dimension $2$
(\cite{nsw}, Cor.~7.2.5).
Assume $U \neq X$ and $i \geq 3$. Then, by \cite{MilADT}, 
Prop~II.2.9., we have $H^i(U,\widehat {\mathcal T})=H^i(G_S,\widehat T)$,
which is zero: indeed $G_S$ is of strict cohomological dimension 
$2$ by \cite{nsw}, Th. 8.3.17.
It remains to deal with the case $U=X$ (now we assume $i \geq 4$). 
By \cite{MilADT}, Lemma~II.2.10, the group
$H^i(X,\Z)$ is torsion; by a restriction-corestriction argument, the 
same holds for $H^i(X,\widehat {\mathcal T})$. Since $\Q$ is uniquely 
divisible, this yields 
$$H^i(X,\widehat {\mathcal T})=H^{i-1}(X,\widehat {\mathcal T} \otimes \Q/\Z)=
\varinjlim_n H^{i-1}(X,\widehat {\mathcal T}/n).$$
For $i \geq 5$, the group $H^{i-1}(X,\widehat {\mathcal T}/n)$ is zero 
(cf. Remark~\ref{fini2}), so we are done. Assume $i=4$.
We observe that the finite group
$H^3(X,\widehat {\mathcal T}/n)$ is dual to $H^0(X,\, _n \T)$ by 
Artin-Mazur-Milne duality (\cite{DemH}, Th.~1.1), so the dual of 
the discrete torsion group $H^4(X,\widehat {\mathcal T})$ is the 
profinite group
$$\varprojlim_n H^0(X,\, _n \T)=\varprojlim_n {_n (T(K))}$$
(the equality holds because the $X$-group scheme $_n \T$ is finite and 
$X$ is connected). But $K$ is a global field, 
hence $T(K)_{\rm tors}$ is finite: indeed if $L$ is a finite extension 
of $K$ such that $T$ splits over $L$, then $T(K) \subset T(L)$ with 
$T(L) \simeq (L^*)^r$ for some $r$, and $L^*$ contains only finitely
many roots of unity. 
Therefore $T(K)$ has trivial Tate module, which yields the result. 

\smallskip

e) Using exact sequence (\ref{dualkummer}), we get a surjection 
$H^2(X,\widehat \T /n) \to \, _n H^3(X, \widehat \T)$, so it is sufficient 
(by Remark~\ref{fini2}) to show that $H^3(X, \widehat \T)$ is of finite
exponent. By restriction-corestriction, we can therefore assume that 
$\widehat \T=\Z$.
By the same method as in d), we get that the dual of 
$H^3(X,\Z)$ is $\varprojlim_n H^1(X,\mu_n)$. As $H^0(X,\G)=k^*$ because 
$X$ is a proper and geometrically integral curve,
we get an exact sequence of finite groups
$$0 \to k^*/k^{*^n} \to H^1(X,\mu_n) \to \, _n \pic X \to 0.$$
Since $\pic X$ is of finite type, we have $\varprojlim_n (_n \pic X)=0$,
hence $\varprojlim_n H^1(X,\mu_n)$ is the inverse limit of the $k^*/k^{*^n}$,
which is $k^*$ itself because $k$ is finite. Thus $H^3(X,\Z)$ is
the dual of $k^*$, which is finite (but not zero).

\enddem

\begin{rema} \label{xcase}
{\rm Assume $U=X$. Then the group $H^2_c(U,\widehat \T)=H^2(X, \widehat \T)$
is in general infinite: for example $H^2(X,\Z) \cong H^1(X,\Q/\Z)$ is the 
dual of the \'etale fundamental group $\pi_1^{\et}(X)$ and the latter is 
an extension of $\gal(\kbar/k)=\widehat \Z$; therefore 
$H^2(X,\Z)$ contains a copy of $\Q/\Z$. 
}
\end{rema}

\begin{prop} \label{finityp}

Let $\C=[\T_1 \stackrel{\rho}{\to} \T_2]$ be a complex of $U$-tori with 
generic fibre $C=[T_1 \to T_2]$.

\smallskip

a) Let $i \in \{-1 , 0 \}$. Then the groups $H^i(U,\widehat {\mathcal C})$
and $H^i_c(U,\widehat {\mathcal C})$ are of finite type, and the 
restriction map $H^i(U,\widehat \C) \to H^i(K,\widehat C)$ is an isomorphism.
The restriction map $H^1(U,\widehat \C) \to H^1(K,\widehat C)$ is injective.
If $U \neq X$, then $H^1_c(U,\widehat {\mathcal C})$ is of finite type. 

\smallskip

b) The groups $H^{-1}(U,\mathcal C)$ and $H^{-1}_c(U,\mathcal C)$ 
are of finite type, and so is $H^0(U,\C)$. 
If $U=X$, then $H^1(U,\C)=H^1(X,\C)$ is of finite type.

\smallskip

c) Assume $U \neq X$. Then $D^1(U,\C)$ and $D^1(U,\widehat \C)$ are finite.

\end{prop}

\dem{} a) The fact that $H^i(U,\widehat {\mathcal C})$
and $H^i_c(U,\widehat {\mathcal C})$ are of finite type for $i \in \{-1,0 \}$
follows immediately by d\'evissage (cf. exact sequences (\ref{devis2})
and (\ref{suppcomp})) from Lemma~\ref{h1lem} a), and we even have 
$H^{-1}_c(U,\widehat {\mathcal C})=0$ if $U \neq X$. For a $U$-torus
$\T$, the restriction map $H^0(U, \widehat \T) \to H^0(K,\widehat T)$ 
obviously is an isomorphism. Let $V$ be a connected Galois covering 
of $U$ with function field $L$ and group $G$, 
such that $\T$ splits over $V$. As seen before,
we have $H^1(V,\Z)=H^1(L,\Z)=0$, hence
$H^1(V,\widehat \T)=H^1(L,\widehat T)=0$. By Hochschild-Serre spectral sequence 
we get $H^1(U,\widehat \T) \cong H^1(K,\widehat T)$ because both groups 
identify to $H^1(G,\widehat T)$. By \cite{MilADT}, Lemma~II.2.10, 
we have 
$$H^2(U,\widehat \T) \cong H^1(U,\widehat \T \otimes \Q/\Z)=
\varinjlim_{n >0} H^1(U, \widehat \T/n),$$
and $H^1(U, \widehat \T/n) \hookrightarrow H^1(K,\widehat T/n)$
because $\widehat \T/n$ is a finite $U$-group scheme, which implies 
$H^2(U,\widehat \T) \hookrightarrow H^2(K,\widehat T)$.

\smallskip 

The commutative diagram with exact lines
{\small 
$$
\begin{CD} 
0 @>>> H^i(U,\widehat \T_1) @>>> H^i(U,\widehat \T_2) @>>>
H^i(U, \widehat \C) @>>>  
H^{i+1}(U,\widehat \T_1) @>>> H^{i+1}(U,\widehat \T_2) \cr
&& @VVV @VVV @VVV @VVV @VVV \cr  
0 @>>> H^i(K,\widehat \T_1) @>>> H^i(K,\widehat T_2) @>>> H^i(K,\widehat C)
@>>> H^{i+1}(K,\widehat T_1) @>>> H^{i+1}(K,\widehat T_2)
\end{CD}
$$
}
and the five lemma
now give that the restriction map $H^{i}(U, \widehat \C) \to
H^{i}(K,\widehat C)$ is an
isomorphism for $i \in \{-1,0 \}$, and is injective for $i=1$.
The fact that $H^1_c(U, \widehat \C)$ is of finite 
type if $U \neq X$ is immediate by d\'evissage thanks to Lemma~\ref{h1lem} c).

\smallskip

b) The first two assertions follow 
from Lemma~\ref{fini1}, using exact sequence (\ref{devis1}). 
For $U=X$, every $X$-torus $T$ satisfies that $H^1(X,\T)$ is of finite type
(Lemma~\ref{fini1}) and $H^2(X,\T)$ is finite (Lemma~\ref{h1lem}, b), 
whence the result. 

\smallskip

c) By functoriality, the image of $D^1(U,{\mathcal C})$ by the
map $u: H^1(U,{\mathcal C}) \to H^2(U,{\mathcal T}_1)$ is a subgroup
of $D^2(U,{\mathcal T}_1)$. The latter is finite by Lemma~\ref{h1lem} b),
because it is a quotient of $H^2_c(U,{\mathcal T}_1)$. As
the kernel of $u$ is a quotient of $H^1(U,\mathcal T_2)$ (which is
finite by Lemma~\ref{fini1} a),
this means that $D^1(U,{\mathcal C})$ is finite.

\smallskip

The group $H^2_c(U,\widehat {\mathcal T}_2)$
is finite by Lemma~\ref{h1lem}, c).
Hence $D^2(U, \widehat {\mathcal T}_2)$ is finite.
On the other hand, the kernel of the map $H^1(U,\widehat {\mathcal C})
\to H^2(U,\widehat {\mathcal T}_2)$ is a quotient of the group
$H^1(U,\widehat {\mathcal T}_1)$, which is finite by Lemma~\ref{h1lem} a).
Thus $D^1(U,\widehat {\mathcal C})$ is finite.

\enddem

\begin{rema} \label{localisom}
{\rm
The same argument as in Proposition~\ref{finityp} 
a) shows that for $v \in U$, the restriction 
map $H^i(\calo_v,\widehat \C)  \to H^i(K_v,\widehat C)$
is an isomorphism for $i \in \{-1, 0 \}$, and is injective for $i=1$. 
}
\end{rema}

Recall (\cite{demphd}, section 3)
that for $v \in X^{(1)}$, given a complex $C$ of $K_v$-tori, the groups $H^{-1}(K_v,C)$ and
$H^0(K_v,C)$ are equipped with a natural Hausdorff topology (and
the groups $H^i(K_v,C)$ are endowed
with the discrete topology for $i \geq 1$, as
are all groups $H^r(K_v, \widehat C)$ for $-1 \leq r \leq 2$).

\begin{lem} \label{closedimage}
The image $I$ of the group $H^0(U,\C)$ into $\bigoplus_{v \not \in U}
H^0(K_v,C)$ is a discrete (hence closed) subgroup, and so is the 
image of $H^{-1}(U,\C)$ into $\bigoplus_{v \not \in U} H^{-1}(K_v,C)$.
\end{lem}

\dem{} We can assume $U \neq X$. Let us start with the case 
when $\C=\G$. 
Then $\calo_U ^*:=H^0(U,\G)$ is a discrete subgroup of $\bigoplus
_{v \not \in U} K_v^*$, because its intersection with the open subgroup
$\bigoplus_{v \not \in U} \calo_v^*$ is $H^0(X,\G)=k^*$, which is finite.
Consider now a $U$-torus $\T$. Let $W$ be a connected Galois 
finite covering of $U$ (with function field $L \supset K$)
that splits $\T$. Let $G:=\gal(L/K)$. 
By the case $\T=\G$, the subgroup $H^0(W,\T)$ is discrete 
in $\bigoplus_{w \not \in W} H^0(L_w,T)$, so $H^0(U,\T)$ is discrete 
in $\bigoplus_{v \not \in U} H^0(K_v,T)$ because it is the intersection 
of $H^0(W,\T) \subset \bigoplus_{w \not \in W} H^0(L_w,T)$ with 
$\bigoplus_{v \not \in U} H^0(K_v,T)=(\bigoplus_{w \not \in W} H^0(L_w,T))^G$.
Thus $I$ is discrete when $\C=\T$ is one single torus.

\smallskip

In the general case, exact triangle (\ref{devistm}) yields a commutative
diagram with exact lines
$$
\begin{CD}
H^1(U,\mathcal M) @>>> H^0(U,\C) @>u>> H^0(U,\T) \cr
@VVV @VVjV @VVV \cr
\bigoplus_{v \not \in U} H^1(K_v,M) @>>> \bigoplus_{v \not \in U}
H^0(K_v,C) @>>> \bigoplus_{v \not \in U} H^0(K_v,T),
\end{CD}
$$
where $\mathcal M$ is a $U$-group of multiplicative type and $\T$
is a $U$-torus. Since $U \neq X$, the right vertical map
is injective. As the lemma holds for $\C=\T$, 
the image $J$ of $H^0(U,\T)$ into
$\bigoplus_{v \not \in U} H^0(K_v,T)$ is discrete,
hence there is an open subgroup $H$ of $\bigoplus_{v \not \in U} H^0(K_v,T)$
such that $J \cap H=\{ 0 \}$. Let $H_1$ be the inverse image of $H$ in
$\bigoplus_{v \not \in U} H^0(K_v,C)$, it is an open subgroup of
$\bigoplus_{v \not \in U} H^0(K_v,C)$ such
that $j^{-1}(H_1)$ is a subgroup of $\ker u$. As
$H^1(U,\mathcal M)$ is finite (Remark~\ref{fini2}), we also have that
$\ker u$ is finite and so is $j^{-1}(H_1)$. Therefore $I \cap H_1
=j(j^{-1}(H_1))$ is finite, which implies that $I$ is discrete.

\smallskip

The same result for the image of $H^{-1}(U,\C)$ into $\bigoplus_{v \not \in U}
H^{-1}(K_v,C)$ follows immediately because $H^{-1}(U,\C)$ is a subgroup 
of $H^0(U,\T_1)$ (which has just been shown to be a discrete subgroup 
of $\bigoplus_{v \not \in U} H^0(K_v,T_1)$), and $\bigoplus_{v \not \in U}
H^{-1}(K_v,C)$ is a topological subspace of
$\bigoplus_{v \not \in U} H^0(K_v,T_1)$.

\enddem

\begin{rema} \label{nocdn}
{\rm The analogue of Lemma~\ref{closedimage} does not hold over a number field
as soon as at least one finite place of $K$ is not in $U$ and $K$ has at
least two non archimedean places: indeed for the exact sequence
$$0 \to H^0(U,\G) \to \bigoplus_{v \not \in U} H^0(K_v,\G) \to
H^1_c(U,\G) \to H^1(U,\G) \to 0$$
to hold (cf. \cite{DemH}, beginning of section 2),
the groups $H^0(K_v,\G)$ at the archimedean places
must be understood as the {\it modified} Tate group ${\widehat H}^0(K_v,\G)$.
The intersection $I$ of $H^0(U,\G)$ with
the compact subgroup 
$\bigoplus_{v \not \in U} \calo_v^* \subset
\bigoplus_{v \not \in U} H^0(K_v,\G) $ (where by convention $\calo_v^*$
means ${\widehat H}^0(K_v,\G)$ at the archimedean places) is countable
and infinite (it is isomorphic to $\calo_K^*$),
hence $I$ is not compact by Baire's Theorem.
Therefore the image of $H^0(U,\G)$ in $\bigoplus_{v \not \in U} H^0(K_v,\G)$
is not closed.
}
\end{rema}

Equip the finitely generated group $H^0(U,\C)$ with the discrete 
topology. We give $H^0_c(U, \C)$ the unique topology such that all maps in 
the exact sequence
\begin{equation} \label {toph0c}
H^{-1}(U,\C)
\to \bigoplus_{v \not \in U} H^{-1}(K_v,C) \to H^0_c(U, \C) \to H^0(U,\C) 
\end{equation} 
are strict (by Lemma~\ref{closedimage}, the left map is strict and the 
quotient of $\bigoplus_{v \not \in U} H^{-1}(K_v,C)$ by the image of 
$H^{-1}(U,\C)$ is a locally compact Hausdorff group).
We also give the finite group (cf. Proposition~\ref{finityp} c)
$D^1(U,\C)$ the discrete topology, and topologize $H^1_c(U,\C)$ such that 
all maps in the exact sequence
\begin{equation} \label {toph1c}
H^0(U,\C) \to \bigoplus_{v \not \in U} H^0(K_v,C) \to H^1_c(U, \C) \to D^1(U,\C) \to 0
\end{equation}
are strict.  

\begin{defi}
{\rm Define ${\mathcal E}$ as the class of those abelian topological groups
$A$ that are an extension
\begin{equation} \label{extgroup}
0 \to P \to A \stackrel{\pi}{\to} F \to 0
\end{equation}
(the maps being continuous)
of a finitely generated group $F$ (equipped with the
discrete topology) by a profinite group $P$ (this implies that all maps
in this exact sequence are strict by \cite{DemH}, Lemma~3.4).
}
\end{defi}

It is easy to check that 
for every group $A$ in $\mathcal E$, a closed subgroup of $A$ and
the quotient of $A$ by any closed subgroup of $A$ are still in $\mathcal E$.
Also a topological extension of a (discrete) finitely generated group by $A$
stays in $\mathcal E$. Finally, every group $A$ in $\mathcal E$ is isomorphic
to the direct product of a finitely generated
group (equipped with discrete topology) by a profinite group: indeed 
up to replacing $F$ by
$F/F_{\rm tors}$ and $P$ by $\pi^{-1}(F_{\rm tors})$
in the extension (\ref{extgroup}), we can assume
that $F=\Z^r$ for some $r \geq 0$. Since $F$ is free, the morphism
$\pi$ has a section $s : F \to A$, which is automatically continuous
because $F$ is discrete. Setting $B=s(F)$, we get a topological isomorphism
$A \cong P \times B$.

\begin{defi}
{\rm
Let $A$ be an abelian group. For every prime number $\ell$,
the {\it $\ell$-adic completion} of $A$ is
$$A^{(\ell)}:=\varprojlim_{m \in \NN^*} (A/\ell^m).$$ We also
set $$A_{\wedge}:=\varprojlim_{n \in \NN^*}
A/n=\prod_{\ell \, {\rm prime}} A^{(\ell)}.$$
The piece of notation $A \{ \ell \}$ stands for the $\ell$-primary
torsion subgroup of $A$.
}
\end{defi}

For $A$ finitely generated, we have
$A^{(\ell)}=A \otimes_{\Z} \Z_{\ell}$; since $\Z_{\ell}$ is a torsion-free
(hence flat) $\Z$-module, the functors $A \mapsto A^{(\ell)}$ and
$A \mapsto A_{\wedge}$ are exact in the category of finitely generated
abelian groups.

\begin{lem} \label{exactcomp}

Let $A \to B \to E \to 0$ be an exact sequence of abelian groups. 

\smallskip

a) The induced map $B^{(\ell)} \to E^{(\ell)}$ is surjective. 

\smallskip

b) If $_{\ell} E$ is finite, then the induced sequence 
$$A^{(\ell)} \to B^{(\ell)} \to E^{(\ell)} \to 0$$
is exact. Likewise if $A/\ell$ is finite.

\end{lem}

\dem{} a) Since the functor $. \otimes_{\Z} \Z/\ell^m$ is right exact, 
the sequence 
$$A/\ell^m \to B/\ell^m \to E/\ell^m \to 0$$
is exact. Therefore the projective system
$(\ker [B/\ell^m \to E/\ell^m])_{m \geq 1}$
has surjective transition maps, which 
implies that the map $\varprojlim_m (B/\ell^m) \to \varprojlim_m (E/\ell^m)$
remains surjective.

\smallskip

b) Assume $E/\ell$ finite. Then 
(by induction on $m$) we also have that $_{\ell^m} E$ is finite
for every positive integer $m$ thanks to the exact sequence 
$$ _{\ell} E \to _{\ell^{m+1}} E \stackrel{.\ell}{\to} _{\ell^m} E.$$
Let $I \subset B$ be the image of $A$ by the map $A \to B$. 
By a), the map $A^{(\ell)} \to I^{(\ell)}$ is surjective, so 
it is sufficient to prove that the sequence 
$$I^{(\ell)} \to B^{(\ell)} \to E^{(\ell)}$$
is exact. By the snake lemma, we have an exact sequence 
$$_{\ell^m} E \to I/\ell^m \to B/\ell^m \to E/\ell^m .$$
Taking projective limit over $m$ yields the required exact sequence 
because the kernel of the map $I/\ell^m \to B/\ell^m$ is finite 
(it is a quotient of $_{\ell^m} E$). Similarly, if $A/\ell$ is finite, 
then $A/\ell^m$ (hence also $I/\ell^m$) is finite for every positive 
$m$ and the same argument works.

\enddem

If we assume further that $A$ is a topological
abelian group, its {\it profinite completion} is
$A^{\wedge}:=\varprojlim_H (A/H)$,
where $H$ runs over all open subgroups of finite index
in $A$. If $A$ is profinite, then $A=A_{\wedge}=A^{\wedge}$. If $A$ is in
the class $\mathcal E$, then $A \hookrightarrow A_{\wedge}=A^{\wedge}$.

\begin{prop} \label{injectcompl}
Let $\mathcal C=[\T_1 \to \T_2]$ be a complex of $U$-tori with generic fiber 
$C=[T_1 \to T_2]$. Let $v \in X^{(1)}$.
The topological groups $H^{-1}(K_v,C)$ and $H^0(K_v,C)$ are 
in $\mathcal E$, as are 
the groups $H^0_c(U,\C)$ and $H^1_c(U,\C)$. In particular, for 
$i \in \{-1, 0 \}$, we have $H^i(K_v,C)_{\wedge}=H^i(K_v,C)^{\wedge}$ 
and for $i \in \{0, 1 \}$, we have $H^i_c(U,\C) \hookrightarrow 
H^i_c(U,\C)_{\wedge}$.
\end{prop}

\dem{}  Let $\T$ be a 
$U$-torus with generic fibre $T$.
Let $v$ be a closed point of $X$ and let $L$ be a finite 
Galois extension of $K_v$ such 
that $\T$ splits over $L$. As $L^* \simeq \Z \times \calo_L^*$
is in $\mathcal E$, so is $T(L)$. 
Then $H^0(K_v,T)$ is in $\mathcal E$
as a closed subgroup of $T(L)$ (the subgroup 
of $\gal(L/K_v)$-invariants). The exact sequence 
$$ H^0(K_v,T_1)  \to H^0(K_v,T_2) \to H^0(K_v,C) \to H^1(K_v,T_1)$$
and the definition of the topology on $H^0(K_v,C)$ now imply that 
$H^0(K_v,C)$ is in $\mathcal E$. The exact sequence (\ref{toph1c}) 
and Lemma~\ref{closedimage} yield that $H^1_c(U,\C)$ is in $\mathcal E$
because $D^1(U,\C)$ is finite.

\smallskip 

Similarly the group $H^{-1}(K_v,C)=\ker [H^0(K_v,T_1) \to H^0(K_v,T_2)]$ is
a closed subgroup of $H^0(K_v,T_1)$, hence is in $\mathcal E$. 
This implies that $H^0_c(U,\C)$ is in $\mathcal E$ thanks to the exact 
sequence (\ref{toph0c}), 
the group $H^0(U,\C)$ being of finite type by Proposition~\ref{finityp}.

\enddem

\begin{lem} \label{vanishlocal}
Let $v \in U$ and let $\mathcal{C}$ be a complex of $\calo_v$-tori. 
Then $H^i(\calo_v, {\mathcal C})=0$ for $i \geq 1$.
\end{lem}

\dem{} Using the exact sequence (\ref{devis1}) with $U$ replaced
by $\calo_v$, we can assume that $\mathcal C=\mathcal T$ is one single
torus. Since $\mathcal T$ is smooth over $\calo_v$, the fppf
cohomology group $H^i(\calo_v, {\mathcal T})$ coincides with the \'etale
group, and it is isomorphic (\cite{MilADT}, Prop.~II.1.1. b) 
to the Galois cohomology group
$H^i(\F_v,\tilde T)$, where $\F_v$
is the residue field of $\calo_v$ and $\tilde T$ the reduction of
$\mathcal T$ mod. $v$, which is a torus over the finite field
$\F_v$. Now $H^1(\F_v,\tilde T)=0$ by Lang's theorem (\cite{lang}).
For $i \geq 2$, the Galois cohomology group $H^i(\F_v,\tilde T)$ is torsion.
Let $n >0$. By Kummer sequence applied to the torus $\tilde T$
over the perfect field $\F_v$, the $n$-torsion
subgroup $_n H^i(\F_v,\tilde T)$ is a quotient of $H^i(\F_v, \, _n T)$, which
is zero because $\F_v$ is of cohomological dimension $1$ (\cite{col},
Chap. XIII, Prop. 2). This proves the lemma.

\enddem

\begin{rema} \label{cdn1}
{\rm Using the definition of fppf compact support cohomology
given in \cite{DemH} (which, in particular, takes care of the 
set $\Omega_{\RR}$ of real places; see loc. cit., Prop.~2.1),
most results of this section hold (with the same proof) if
we replace $K$ by a number field with ring of integers $\calo_K$,
$X$ by $\spec \calo_K$, and
$U$ by a non empty Zariski open subset of $X$. Also the piece of notation 
$v \not \in U$ means that we consider the closed points of $X-U$ {\it and} 
the real places of $K$; for $v \in \Omega_{\RR}$ and $i \leq 0$,
the groups $H^i(K_v,...)$ 
must be understood as the modified Tate groups (cf. Remark~\ref{nocdn}).
More precisely:

\begin{itemize}

\item Lemma~\ref{fini1} a) hold without the restriction $U \neq X$; b) 
and Remark~\ref{fini2} are useless because for every finite flat commutative 
$U$-group scheme $\mathcal N$, all groups $H^i(U,\mathcal N)$
and $H^i_c(U,\mathcal N)$ are finite (cf. \cite{DemH}, Th.1.1). 

\item Lemma~\ref{h1lem} a) and b) are unchanged; c) holds without the 
restriction $U \neq X$. In d), the vanishing of $H^i_c(U,\widehat \T)$ 
for $i \geq 4$ still holds but the proof uses a different argument (namely 
that the dual of this group is the inverse limit of the $H^{4-i}(U, \, _n \T)$,
which is zero even in the case $i=4$ because the finitely generated 
group $H^0(U,\T)$ has trivial Tate module). The vanishing
of $H^i(U,\widehat \T)$ for $i \geq 4$ must be replaced by its finiteness if  
$\Omega_{\RR} \neq \emptyset$. 
Finally, {\it the vanishing of $H^3(U,\widehat \T)$ does not 
hold anymore}, even if $\Omega_{\RR}=\emptyset$
(if Leopoldt's conjecture is assumed, then 
$H^3(U,\widehat \T) \{ \ell \}=0$ for $\ell$ invertible on $U$, but
not in general for other primes; \cite{MilADT}, Th. II.4.6. b) is wrong for 
$r=3$, the problem in the proof being that the second line of the diagram 
needs not remain exact after taking profinite completions); neither does e)
hold as soon as there are at least two archimedean places.
Also, there is no more counterexample as in Remark~\ref{xcase}. 

\item Proposition~\ref{finityp} a) and c) 
hold (without the condition $U \neq X$), 
b) is also true (and when $U=X$, the group $H^1(X,\C)$ is even finite).

\end{itemize}

}
\end{rema}

\section{Duality theorems in fppf cohomology} \label{two}

In order to state and prove duality results for the cohomology of complexes of fppf sheaves, we need to extend some constructions from \cite{DemH} to the context of bounded complexes.
Let $A$ and $B$ be two bounded complexes of fppf sheaves over $U$. Following \cite{SGA4}, XVII, 4.2.9 or \cite{friedsus}, appendix A, one can consider the Godement resolutions $G(A)$ and $G(B)$ of $A$ and $B$. As in \cite{God}, II.6.6 or in \cite{friedsus}, appendix A, there is a natural commutative diagram of complexes
\[
\xymatrix{
A \otimes B \ar[rd] \ar[d] & \\
\textup{Tot}(G(A) \otimes G(B)) \ar[r] & G(A \otimes B) \, .
}
\]
Following \cite{DemH}, proof of Lemma 4.1, one gets a functorial morphism of complexes
\begin{equation} \label{cup complexes} 
\textup{Tot}(\Gamma_c(U, G(A)) \otimes \Gamma(U,G(B))) \to \Gamma_c(U,G(A \otimes B))
\end{equation}
and a functorial pairing
$$R\Gamma_c(U,A) \otimes^{\L} R\Gamma(U,B) \to R\Gamma_c(U,A \otimes B) \, .$$

In particular, if $C$ is a bounded complex and $\widehat{C}  := \underline{\textup{Hom}}^\bullet(C, \Gm[1])$ its dual, using the morphism $\textup{Tot}(C \otimes \widehat{C}) \to \Gm[1]$ from section \ref{one}, we get functorial pairings
\[\textup{Tot}(\Gamma_c(U, G(C)) \otimes \Gamma(U,G(\widehat{C}))) \to \Gamma_c(U,G(\Gm[1]))\]
and
\begin{equation} \label{pairing derived}
R\Gamma_c(U,C) \otimes^{\L} R\Gamma(U,\widehat{C}) \to R\Gamma_c(U,\Gm[1]) \, .
\end{equation}

\smallskip

Following section \ref{one}, given a bounded complex $\mathcal{C}$ of finite flat commutative group schemes over $U$, there is a natural topology on the abelian groups $H^i_c(U, \mathcal{C})$. This topology is profinite via Lemma \ref{lem topo profinite}, and considering $H^j(U, \widehat{\mathcal{C}})$ as discrete torsion groups, the pairings
\[
H^i_c(U, \mathcal{C}) \times H^j(U, \widehat{\mathcal{C}}) \to H^{i+j+1}_c(U, \G)
\]
are continuous by the same argument as \cite{DemH}, Lemma~4.4.

\begin{prop} \label{AMM complex}
Let $\mathcal{C}$ be a bounded complex of finite flat commutative group schemes over $U$.
For all $i \in \Z$, there is a perfect pairing between profinite and discrete torsion groups 
\[
H^{2-i}_c(U, \mathcal{C}) \times H^i(U, \widehat{\mathcal{C}}) \to H^3_c(U,\G) \cong \Q/\Z \, .\]
\end{prop}

\dem{}
The isomorphism $H^3_c(U,\G) \cong \Q/\Z$ follows from \cite{MilADT}, \S 2.3 and \cite{DemH}, Prop.~2.1., 4).

We now prove the Proposition by induction on the length of the complex $\mathcal{C}$.
\begin{itemize}
	\item if $\mathcal{C}$ is concentrated in a given degree $n$, then the Proposition is a direct consequence of \cite{DemH}, Theorem 1.1.
	\item assume that $\mathcal{C} := [C_r \xrightarrow{f_r} C_{r+11} \to \dots \to C_s]$, with $C_i$ in degree $i$, has length $s-r \geq 1$.
Let $A := \ker(f_r)$, and let $\overline{\mathcal{C}} := \textup{Cone}(A[-r] \xrightarrow{j} \mathcal{C})$. Then there is an exact triangle
\[
A[-r] \xrightarrow{j} \mathcal{C} \xrightarrow{i} \overline{\mathcal{C}} \xrightarrow{p} A[1-r] \, .
\]
We apply the functor $\widehat{\cdot}$ to this exact triangle. Then, using \eqref{cone and dual}, we get that the natural triangle
\[
\widehat{A[1-r]} \xrightarrow{\widehat{p}} \widehat{\overline{\mathcal{C}}} \xrightarrow{\widehat{i}} \widehat{\mathcal{C}} \xrightarrow{\widehat{j}} \widehat{A[-r]}  
\]
is exact. 

Since the pairing between a complex and its dual is functorial, we get from the previous triangles a commutative diagram of topological groups with exact rows, where the vertical maps comes from the pairings \eqref{pairing derived}:
\begin{equation} \label{devissage AMM1}
\xymatrix{
H^{i-1}_c(U, \overline{\mathcal{C}}) \ar[r] \ar[d] & H^{i-r}_c(U, A) \ar[r] \ar[d] & H^i_c(U, \mathcal{C}) \ar[r] \ar[d] & H^{i}_c(U, \overline{\mathcal{C}}) \ar[r] \ar[d] & H^{i+1-r}_c(U, A) \ar[d] \\
H^{3-i}(U, \widehat{\overline{\mathcal{C}}})^* \ar[r] & H^{3+r-i}(U, A^D)^* \ar[r] & H^{2-i}(U, \widehat{\mathcal{C}})^* \ar[r] & H^{2-i}(U, \widehat{\overline{\mathcal{C}}})^* \ar[r] & H^{2+r-i}(U, A^D)^* \, . \\
}
\end{equation}
The second line remains exact as the dual of an exact sequence of discrete groups.

Note that we have a quasi-isomorphism $\varphi : \overline{\mathcal{C}} \to \mathcal{C}'$, where $\mathcal{C}' := [C_{r+1}/Im(f_r) \to C_{r+2} \to \dots \to C_s]$ has a smaller length than $\mathcal{C}$, hence by induction, we can assume that the natural maps $H^r_c(U, \mathcal{C}') \to H^{2-r}(U, \widehat{\mathcal{C}'})^*$ are isomorphisms. Since all the $C_i$'s are finite flat group schemes, we see that the dual morphism $\widehat{\varphi} : \widehat{\mathcal{C}'} \to \widehat{\overline{\mathcal{C}}}$ is a quasi-isomorphism, hence by functoriality of the pairings, we deduce that the maps $H^r_c(U, \overline{\mathcal{C}}) \to H^{2-r}(U, \widehat{\overline{\mathcal{C}}})^*$ are isomorphisms too.

Hence in diagram \eqref{devissage AMM1}, all vertical morphisms, except perhaps the central one, are isomorphisms. Then the five Lemma implies that the central morphism is an isomorphism.
\end{itemize}
By induction on the length of $\mathcal{C}$, the Proposition is proven.
\enddem

{\it From now on, we denote by $\C=[\T_1 \to \T_2]$ a complex of $U$-tori} 
with generic fibre $C=[T_1 \to T_2]$ and dual $\widehat {\mathcal C}$ (cf. 
section \ref{onebis}).
As a consequence of the previous proposition, we can get the global duality results for the cohomology of the complex $\mathcal{C} \otimes^{\bf L} \Z /n$: 

\begin{prop} \label{avfini}
Let $n$ be a positive integer (not necessarily prime to $p$).
Let $i$ be an integer with $-2 \leq i \leq 2$. 
\begin{enumerate}
	\item There is a perfect pairing of finite groups 
$$H^i(U, {\mathcal C} \otimes^{\bf L} \Z/n) \times H^{1-i}_c(U, 
\widehat {\mathcal C} \otimes^{\bf L} \Z/n) \to \Q/\Z \, .$$

If $U \neq X$, the groups $H^2(U, {\mathcal C} \otimes^{\bf L} \Z/n)$ and 
$H^{-1}_c(U, \widehat {\mathcal C} \otimes^{\bf L} \Z/n)$ are zero.
	\item There is a perfect pairing  
$$H^i_c(U, {\mathcal C} \otimes^{\bf L} \Z/n) \times H^{1-i}(U, 
\widehat {\mathcal C} \otimes^{\bf L} \Z/n) \to \Q/\Z$$
between the profinite group $H^i_c(U, {\mathcal C} \otimes^{\bf L} \Z/n)$ 
and the discrete group $H^{1-i}(U, \widehat {\mathcal C} \otimes^{\bf L} \Z/n)$. 
These groups are finite if $i \not \in \{ 0, 1 \}$ or if $p$ and $n$ are 
coprime. The groups
$H^{-2}_c(U,{\mathcal C} \otimes^{\bf L} \Z/n) $ and 
$H^3(U,\widehat {\mathcal C} \otimes^{\bf L} \Z/n)$ are zero if 
$U \neq X$.
\end{enumerate} 
Moreover, all the groups involved are zero if $\vert i \vert >2$.
\end{prop}

\dem{} 
Recall that there is a quasi-isomorphism of complexes $\psi = C' \to \mathcal{C} \otimes^{\bf L} \Z / n$, where $C' := [{_n T_1} \xrightarrow{\rho} {_n T_2}]$, with ${_n T_1}$ in degree $-2$, and that the dual morphism $\widehat{\psi} : \widehat{\mathcal{C} \otimes^{\bf L} \Z / n} \to \widehat{C'} = [\widehat{T_2}/n \xrightarrow{-\widehat{\rho}} \widehat{T_1}/n] = \left(\widehat{\mathcal{C}} \otimes^{\bf L} \Z / n \right)[-1]$ (with $\widehat{T_2}/n$ in degree $0$) is also a quasi-isomorphism.

Since $C'$ is a bounded complex of finite flat commutative group schemes, then Proposition \ref{AMM complex} implies that the pairings in the statement of the Proposition are perfect pairings of topological groups.

Let us now check the finiteness and vanishing results.
Using the exact triangle (\ref{exactriang1}), we get an exact sequence:
\[
H^{i+2}(U,{_n \ker \rho}) \to H^i(U, {\mathcal C} \otimes^{\bf L} \Z/n) \to H^{i+1}(U, T_{\Z/n}({\mathcal C})) \, .
\]
As the finite $U$-group schemes $T_{\Z/n}({\mathcal C})$ and 
$_n \ker \rho$ are of multiplicative type, Lemma~\ref{fini1} 
and Remark~\ref{fini2}
imply that the groups $H^r(U,T_{\Z/n}({\mathcal C}))$ and 
$H^r(U,\, _n \ker \rho)$ are finite for every integer $r$. Thus 
$H^i(U, {\mathcal C} \otimes^{\bf L} \Z/n)$ is finite.

\smallskip

Similarly, using the exact triangle (\ref{exactriang1}), the group $H^i_c(U, {\mathcal C} \otimes^{\bf L} \Z/n)$ is finite for
$i \not \in \{ 0, 1 \}$ (or if $p$ and $n$ are coprime) by Remark~\ref{fini2}. 

\smallskip

Recall that for a finite and flat group scheme $\mathcal N$ over $U$, 
we have $H^r(U,\mathcal N)=0$ for $r <0$ (obvious), for $r \geq 4$, 
and also for $r=3$ if $U \neq X$~: indeed by loc. cit. $H^r(U,\mathcal N)$
is dual to $H^{3-r}_c(U,{\mathcal N}^D)$; the latter is clearly zero if 
$r \geq 4$ and
if $U \neq X$, we also have $H^0_c(U,{\mathcal N}^D)=0$
thanks to the exact sequence 
$$0 \to H^0_c(U,{\mathcal N}^D) \to H^0(U,{\mathcal N}^D) \to 
\bigoplus_{v \not \in U} H^0(K_v,N^D),$$
the last map being injective by the assumption $U \neq X$. 
The previous d\'evissages now yield the vanishing assertions of 
the proposition.
\enddem

\begin{prop} \label{finivanish1}

Let $i$ be an integer. Let $n$ be a positive integer (not necessarily prime 
to $p$).

\smallskip

a) The groups $_n H^i(U, {\mathcal C})$ and 
$H^i(U, {\mathcal C})/n$ are finite. The group 
$H^i(U, {\mathcal C})$ is torsion if $i \geq 2 $ (resp. if $i \geq 1$ and 
$U \neq X$). Besides $H^i(U, {\mathcal C})=0$ in the following cases~:
$i \geq 4$; $i \leq -2$; $i=3$ and $U \neq X$.  

\smallskip

b) The groups $_n H^i_c(U, \widehat {\mathcal C})$  and 
$H^i_c(U, \widehat {\mathcal C})/n$ are finite. The group 
$H^i_c(U, \widehat {\mathcal C})$ is torsion if $i \geq 2$, and it is 
zero if $i \geq 4$ or $i \leq -2$. Assume further $U \neq X$; then 
$H^i_c(U, \widehat {\mathcal C})=0$ for $i=-1$.

\end{prop}

\dem{} a) Using the exact triangle in ${\mathcal D}^b(U)$~:
$$ {\mathcal C} \stackrel{.n}{\to} {\mathcal C} \to {\mathcal C}
\otimes^{\bf L} \Z/n \to {\mathcal C}[1] ,$$
we get an exact sequence of abelian groups 
\begin{equation} \label {kummer1}
0 \to H^i(U,{\mathcal C})/n \to H^i(U,{\mathcal C} \otimes^{\bf L} \Z/n) 
\to _n H^{i+1}(U,,{\mathcal C}) \to 0.
\end{equation}
since $H^i(U,{\mathcal C} \otimes^{\bf L} \Z/n)$ is finite
(Proposition~\ref{avfini} 1.), the finiteness of the groups
$_n H^{i+1}(U, {\mathcal C})$ and $H^i(U, {\mathcal C})/n$ follows 
for all $i$.

\smallskip

To prove that $H^i(U, {\mathcal C})$ is torsion if $i \geq 2 $ (resp. 
if $U \neq X$ and $i=1$), we can restrict by d\'evissage (using 
exact sequence (\ref{devis1})) to the case when 
${\mathcal C}={\mathcal T}$ is one single torus. If $U \neq X$ and $i=1$,
this follows from Lemma~\ref{fini1}, a), so assume $i \geq 2$. 
We can also assume by a restriction-corestriction argument 
that ${\mathcal T}=\G$
because the torus ${\mathcal T}$ is split by some finite 
\'etale covering of $U$. 
Now $H^2(U,\G)=\br U$ is torsion because it injects into $\br K$; 
also $H^3(U,\G)$ is torsion (it is 
even $0$ if $U \neq X$) and $H^i(U,\G)=0$ for $i \geq 4$ by \cite{MilADT}, 
Prop. II.2.1., the group scheme $\G$ being smooth (hence \'etale and fppf 
cohomology coincide). 

\smallskip For every $U$-torus ${\mathcal T}$, we have 
$H^i(U,{\mathcal T})=0$ for negative $i$ (obvious), hence by d\'evissage
$H^i(U,{\mathcal C})=0$ for $i <-1$. Let $i \geq 3$; as seen before
$H^i(U,{\mathcal C})$ is torsion and $_n H^i(U,{\mathcal C})$ is 
a quotient of $H^{i-1}(U,{\mathcal C} \otimes^{\bf L} \Z/n)$ by exact 
sequence (\ref{kummer1}). The latter is zero if $i \geq 4$, and also
if $i=3$ when $U \neq X$ by 
the vanishing assertions in Proposition~\ref{avfini} 1.. Thus 
$H^i(U,{\mathcal C})$ is zero if $i \geq 4$, and also if $i=3$ if we 
assume further $U \neq X$.

\smallskip
 
b) Similarly, the finiteness statements follow from the exact sequence 
\begin{equation} \label{kummer2}
0 \to H^i_c(U,\widehat {\mathcal C})/n \to H^i_c(U,\widehat {\mathcal C} 
\otimes^{\bf L} \Z/n) \to \, _n H^{i+1}_c(U,\widehat {\mathcal C}) \to 0
\end{equation}
combined to Proposition~\ref{avfini} 1.. Let $i \geq 2$. To prove that 
$H^i_c(U,\widehat {\mathcal C})$ is torsion we can assume that $\widehat
{\mathcal C}$ is the dual of a torus (via exact sequence (\ref{devis2})), 
then that 
$\widehat {\mathcal C}=\Z$ (by a restriction-corestriction argument). 
Using the exact sequence 
$$\bigoplus_{v \not \in U} H^{i-1}(K_v, \Z) \to H^i_c(U, \Z) \to 
H^i(U, \Z),$$ 
it is sufficient to prove that $H^i(U,\Z)$ is torsion because the Galois 
cohomology groups $H^{i-1}(K_v, \Z)$ are torsion for $i-1>0$. This holds
by \cite{MilADT}, Lemma~II.2.10. 

\smallskip Let ${\mathcal T}$ be a $U$-torus. For each integer $i$, 
there is an exact sequence 
\begin{equation} \label{suppc} 
\bigoplus_{v \not \in U} H^{i-1}(K_v, \widehat T) \to 
H^i_c(U, \widehat {\mathcal T}) \to H^i(U, \widehat {\mathcal T}) \to 
\bigoplus _{v \not \in U} H^i(K_v, \widehat T).
\end{equation}
Therefore $H^i_c(U,\widehat {\mathcal T})=0$ for $i <0$, hence 
$H^i_c(U,\widehat {\mathcal C})=0$ (by d\'evissage) for $i <-1$. 
For $i \geq 4$, we have
$H^i_c(U,\widehat {\mathcal T})=0$ by Lemma~\ref{h1lem} d), and 
$H^i_c(U,\widehat {\mathcal C})=0$ by d\'evissage.

\smallskip

Assume now $U \neq X$. Then $H^0_c(U,\widehat {\mathcal T})=0$ by exact sequence
(\ref{suppc}) applied to $i=0$: indeed the map 
$H^0(U,\widehat {\mathcal T}) \to \bigoplus_{v \not \in U} 
H^0(K_v,\widehat T)$ is injective (choose a closed point $v$ of $X-U$; then 
the restriction maps $H^0(U, \widehat {\mathcal T}) \to H^0(K,\widehat T)$
and $H^0(K,\widehat T) \to H^0(K_v,\widehat T)$ are injective).
Therefore $H^{-1}_c(U, \widehat {\mathcal C})=0$
by d\'evissage.  

\enddem

\begin{prop} \label{finivanish2}
Let $i$ be an integer. Let $n$ be a positive integer (not necessarily prime
to $p$).

\smallskip

a) The group $_n H^i(U,\widehat {\mathcal C})$ is finite if $i \not \in 
\{1, 2 \}$. The group $H^i(U,\widehat {\mathcal C})/n$ is finite if
$i \neq 1$. The group $H^i(U,\widehat {\mathcal C})$ is torsion 
if $i \geq 1$, and it is zero if $i \geq 4$ or $i \leq -2$. If we assume 
further $U \neq X$, then $H^3(U,\widehat {\mathcal C})=0$.

\smallskip 

b) The group $_n H^i_c(U,{\mathcal C})$ is finite if $i \neq 1$. 
The group $H^i_c(U,{\mathcal C})/n$ is finite if 
$i \not \in \{ 0,1 \}$.
The group $H^i_c(U,{\mathcal C})$ is torsion for $i \geq 2$, and it is 
zero if $i \geq 4$ or $i \leq -2$ (resp. $i=-1$ if $U \neq X$).

\end{prop}

\dem{} a) The exact sequence 
\begin{equation} \label{kummer3}
0 \to H^i(U,\widehat {\mathcal C})/n \to H^i(U,\widehat {\mathcal C} 
\otimes^{\bf L} \Z/n) \to _n H^{i+1}(U,\widehat {\mathcal C}) \to 0
\end{equation}
and Proposition~\ref{avfini} 2. yield the finiteness of 
$_n H^i(U,\widehat {\mathcal C})$ for $i \not \in \{1, 2 \}$
and of $H^i(U,\widehat {\mathcal C})/n$ for $i \not \in \{0, 1 \}$.
Besides the abelian group $H^0(U,\widehat {\mathcal C})$ is of finite type
by Proposition~\ref{finityp}.  
In particular $H^0(U,\widehat {\mathcal C})/n$ is finite.

\smallskip

Let $i \geq 1$. To prove that $H^i(U,\widehat {\mathcal C})$ is torsion,
we can assume (by d\'evissage) that $\widehat \C=\widehat \T$, where 
$\T$ is one single torus,
then that $\widehat {\mathcal C}=\Z$ (by restriction-corestriction); then 
the result holds by \cite{MilADT}, Lemma~II.2.10. For $i \leq -2$ or 
$i \geq 4$ (resp. $i=3$ if $U \neq X$), 
the group $H^i(U,\widehat {\mathcal C})$ is zero by 
d\'evissage (using Lemma~\ref{h1lem}, d) for the latter).

\smallskip

b) There is an exact sequence 
\begin{equation} \label{kummer4}
0 \to H^i_c(U,{\mathcal C})/n \to H^i_c(U,{\mathcal C} 
\otimes^{\bf L} \Z/n) \to _n H^{i+1}_c(U,{\mathcal C}) \to 0
\end{equation}
By Proposition~\ref{avfini} 2., the group $_n H^i_c(U,{\mathcal C})$ 
is finite for $i \not \in \{1 , 2 \}$ and the group 
$H^i_c(U,{\mathcal C})/n$ is finite for $i \not \in \{0 , 1 \}$.
To prove that the groups 
$H^i_c(U,{\mathcal C})$ are torsion for $i \geq 2$, we can assume 
as usual that ${\mathcal C}={\mathcal T}$ is a torus.
Then we apply 
Proposition~\ref{finivanish1} a) and the exact sequence 
$$\bigoplus_{v \not \in U} H^{i-1}(K_v,T) \to 
H^i_c(U,{\mathcal T}) \to H^i(U,{\mathcal T}) .$$

\smallskip

Besides $H^2_c(U,{\mathcal T})$ is finite  by Lemma~\ref{h1lem} b), so
$H^2_c(U,{\mathcal C})$ is also of cofinite type by d\'evissage 
because we already know that 
$H^3_c(U,{\mathcal T}_1)$ is torsion of cofinite type.

\smallskip

Obviously we have $H^i_c(U,{\mathcal T})=0$ for every negative $i$, hence 
$H^i_c(U,{\mathcal C})=0$ by d\'evissage if $i \leq -2$. Assume $i \geq 4$. 
then $H^i_c(U,{\mathcal C})=H^i(U,{\mathcal C})$ (apply 
exact sequence (\ref{suppcomp}) and use the fact that $K_v$ is of 
strict cohomological dimension 2), so $H^i_c(U,{\mathcal C})=0$
by Proposition~\ref{finivanish1} a). If we assume further $U \neq X$, 
then $H^0_c(U,{\mathcal T})=0$ (same argument as in
Proposition~\ref{finivanish1} b), so $H^{-1}_c(U,{\mathcal C})=0$ by d\'evissage.

\enddem

\begin{rema} \label{nocofini}
{\rm Using Remark~\ref{fini2}, it is easy to see that the finiteness 
assertions of Proposition~\ref{finivanish2} hold for every $i$ if we assume 
further that $p$ does not divide $n$, but this is no longer true 
in general if $U \neq X$.
Indeed the group $H^1(U,\Z/p)$ and its dual $H^2_c(U,\mu_p)$ 
can be infinite (cf. \cite{MilADT}, Lemma~III.8.9). 
Since $H^1(U,\Z)=0$ and $_p H^2_c(U,\G)$ is finite,
this implies that $_p H^2(U,\Z)$ and $H^1_c(U,\G)/p$ are infinite, which gives 
examples of $_p H^i(U, \widehat \C)$ infinite for $i=1,2$ and of 
$H^i_c(U,\C)/p$ infinite for $i=0,1$. The complex
$\C=[\G \stackrel{.p}{\to} \G]$ is an example with $_p H^1_c(U,\C)$ 
and $H^1(U,\widehat \C)/p$ infinite (indeed $\C$ is quasi-isomorphic 
to $\mu_p[1]$ and $\widehat \C$ is quasi-isomorphic to $\Z/p$).
}
\end{rema}

\begin{rema}
{\rm For every integer $r$ and every positive integer $n$, 
the groups $H^r(U, \mathcal C)/n$ and 
$H^r_c(U,\widehat {\mathcal C})/n$ are finite by
Proposition~\ref{finivanish1} b), so for each prime number 
$\ell$ (including $\ell=p$), the $\ell$-adic completions 
$$H^r(U,{\mathcal C})^{(\ell)}:=\varprojlim_m 
H^r(U,{\mathcal C})/\ell^m; \quad 
H^r_c(U,\widehat {\mathcal C})^{(\ell)}:=\varprojlim_m 
H^r_c(U,\widehat {\mathcal C})/\ell^m $$
are profinite. Exact sequence (\ref{kummer4}) 
shows that $H^r_c(U,{\mathcal C})/n$ is a closed subgroup 
of the profinite group $H^r_c(U,{\mathcal C} \otimes^{\bf L} \Z/n)$,
hence $H^r_c(U,{\mathcal C})/n$ is profinite and so is the 
$\ell$-adic completion $H^r_c(U,{\mathcal C})^{(\ell)}$.
}
\end{rema}

The map ${\mathcal C} \otimes^{\bf L} \widehat {\mathcal C} \to \G[1]$
induces for every integer $r$ pairings
\begin{equation} \label{leftpairing}
H^r(U, {\mathcal C}) \times H^{2-r}_c(U,\widehat {\mathcal C}) \to
H^3_c(U,\G) \cong \Q/\Z.
\end{equation}

\begin{equation} \label{rightpairing}
H^r_c(U, {\mathcal C}) \times H^{2-r}(U,\widehat {\mathcal C}) \to
H^3_c(U,\G) \cong \Q/\Z.
\end{equation}

We now prove a key lemma.

\begin{lem} \label{avlem}
Let $\ell$ be a prime number (possibly equal to $p$). Let $i$ be 
an integer. 

\smallskip

a) The maps $$\psi : H^{i+1}(U,{\mathcal C})\{ \ell \} \to
(H^{1-i}_c(U,\widehat {\mathcal C})^{(\ell)})^*$$
$$\psi' : H^{1-i}_c(U,\widehat {\mathcal C})\{ \ell \} \to 
(H^{i+1}(U,{\mathcal C})^{(\ell)})^*$$
induced by the pairing (\ref{leftpairing}) are surjective 
and have divisible kernel. Besides $\psi'$ is an isomorphism 
if $i=-2$ and $\psi$ is an isomorphism if we have 
both $i=1$ and $U \neq X$. 

\smallskip

b)  The map $$\varphi : H^{i+1}(U,\widehat {\mathcal C})\{ \ell \} \to
(H^{1-i}_c(U,{\mathcal C})^{(\ell)})^*$$
induced by the pairing (\ref{rightpairing}) is surjective, and
has divisible kernel (resp. is an isomorphism if we assume 
both $i=1$ and $U=X$). 
Assume $i \not \in \{ -1 , 0 \}$. Then the map 
$$\varphi' : H^{1-i}_c(U,{\mathcal C})\{ \ell \} \to
(H^{i+1}(U,\widehat {\mathcal C})^{(\ell)})^*$$
is surjective and has divisible kernel (resp. is an isomorphism 
if $i=-2$).

\end{lem}

Observe that for $U \neq X$, the groups involved can be non zero only if
$-2 \leq i \leq 1$. 

\dem{} a) For each positive integer $m$, 
there is an exact  commutative diagram of finite abelian groups, 
\begin{equation} \label{finidiag1}
\begin{CD}
0 @>>> H^i(U,{\mathcal C})/\ell^m @>>> H^i(U,{\mathcal C} \otimes^{\bf L}
 \Z/\ell^m) @>>> _{\ell^m} H^{i+1}(U,{\mathcal C}) @>>> 0 \cr 
&& @VVV @VVV @VV{\psi_m}V \cr 
0 @>>>( _{\ell^m} H^{2-i}_c(U,\widehat {\mathcal C}))^* @>>> H^{1-i}_c(U,
\widehat {\mathcal C} \otimes^{\bf L}
 \Z/\ell^m)^* @>>> (H^{1-i}_c(U,\widehat {\mathcal C}) /\ell^m)^* @>>> 0 .
\end{CD}
\end{equation}
By Proposition~\ref{avfini} 1., the middle vertical map is an isomorphism. 
Taking direct limit over $m$ and applying the snake lemma, we 
get that $\psi=\varinjlim_m \psi_m$
is surjective and ${\rm Ker} \, \psi$ is a quotient
of $(T_{\ell} H^{2-i}_c(U,\widehat {\mathcal C}))^* $. Since each
$_{\ell^m} H^{2-i}_c(U,\widehat {\mathcal C})$ is finite, 
the $\ell$-adic Tate module $T_{\ell} H^{2-i}_c(U,\widehat {\mathcal C})$ 
is profinite and torsion-free, which implies that 
the dual $(T_{\ell} H^{2-i}_c(U,\widehat {\mathcal C}))^*$  
is divisible. 
For $U \neq X$ and $i=1$, the group $H^1_c(U,\widehat {\mathcal C})$
is of finite type by Proposition~\ref{finityp} a), so its 
$\ell$-adic Tate-module is zero and $\psi$ has trivial kernel.

\smallskip

The argument for $\psi'$ is similar, using the exact commutative diagram
\begin{equation} \label{finidiag2}
\begin{CD}
0 @>>> H^{-i}_c(U,\widehat {\mathcal C})/\ell^m @>>> H^{-i}_c(U, \widehat 
{\mathcal C} \otimes^{\bf L}
 \Z/\ell^m) @>>> _{\ell^m} H^{1-i}_c(U,\widehat {\mathcal C}) @>>> 0 \cr 
&& @VVV @VVV @VV{\psi'_m}V \cr 
0 @>>>( _{\ell^m} H^{2+i}(U,{\mathcal C}))^* @>>> H^{2+i}(U,
{\mathcal C} \otimes^{\bf L}
 \Z/\ell^m)^* @>>> (H^{1+i}(U,{\mathcal C}) /\ell^m)^* @>>> 0 .
\end{CD}
\end{equation}

Besides, for $i=-2$, the Tate module 
of $H^{2+i}(U,C)=H^0(U,\mathcal C)$ is trivial because 
$H^0(U,\mathcal C)$ is a finitely generated abelian group
(Proposition~\ref{finityp}, b), which gives that $\psi'$ has trivial 
kernel.  

\smallskip

b) There is an exact commutative diagram of discrete 
abelian groups (observe that the second line is obtained by 
dualizing an exact sequence of profinite groups):
\begin{equation} \label{finidiag3}
\begin{CD} 
0 @>>> H^i(U,\widehat {\mathcal C})/\ell^m @>>> H^i(U,\widehat 
{\mathcal C} \otimes^{\bf L}
 \Z/\ell^m) @>>> _{\ell^m} H^{i+1}(U,\widehat {\mathcal C}) @>>> 0 \cr
&& @VVV @VVV @VV{\varphi_m}V \cr
0 @>>>( _{\ell^m} H^{2-i}_c(U,{\mathcal C}))^* @>>> H^{1-i}_c(U,
{\mathcal C} \otimes^{\bf L}
 \Z/\ell^m)^* @>>> (H^{1-i}_c(U,{\mathcal C}) /\ell^m)^* @>>> 0 .
\end{CD}
\end{equation}
Since the middle vertical is an isomorphism by Proposition~\ref{avfini} 2.
and $_{\ell^m} (H^{2-i}_c(U,{\mathcal C}))$ is profinite for each $m$ 
(hence $H^{2-i}_c(U,{\mathcal C})$ has profinite $\ell$-adic Tate 
module), the same argument as in a) yields that $\varphi$ is surjective with 
divisible kernel. If $U=X$ and $i=1$, then $H^{2-i}_c(U,{\mathcal C})=
H^1(X,\C)$ is of finite type by Proposition~\ref{finityp} b),
so it has trivial $\ell$-adic Tate module and $\varphi$ is an isomorphism.

\smallskip

The argument for $\varphi'$ is similar, except that we use the 
exact commutative diagram 
\begin{equation} \label{finidiag4}
\begin{CD} 
0 @>>> H^{-i}_c(U,{\mathcal C})/\ell^m @>>> H^{-i}_c(U, 
{\mathcal C} \otimes^{\bf L}
 \Z/\ell^m) @>>> _{\ell^m} H^{1-i}_c(U,{\mathcal C}) @>>> 0 \cr
&& @VVV @VVV @VV{\varphi'_m}V \cr
0 @>>>( _{\ell^m} H^{i+2}(U,\widehat {\mathcal C}))^* @>>> H^{i+1}(U,
\widehat {\mathcal C} \otimes^{\bf L}
 \Z/\ell^m)^* @>>> (H^{i+1}(U,\widehat {\mathcal C}) /\ell^m)^* @>>> 0 .
\end{CD}
\end{equation}

only for $i \not \in \{ -1 , 0 \}$ (for $i\in \{-1, 0\}$ and $U \neq X$,
the diagram would consist of profinite but possibly infinite groups if 
$\ell=p$, so  
direct limits would not necessarily behave well; in particular, $\ell$-adic 
completions involved would not necessarily be profinite). The same 
argument as in a) shows that $\varphi'$ is surjective with divisible 
kernel, and this kernel is trivial for $i=-2$ because the finitely generated 
abelian group $H^0(U, \widehat {\mathcal C})$ (cf. Proposition~\ref{finityp}, 
a) has trivial $\ell$-adic Tate module.  

\enddem

\begin{rema}
{\rm For abelian groups $A$, $B$, assertions like ``$A \{\ell \} \to 
(B^{(\ell)}) ^*$ is surjective with divisible kernel'' can be rephrased 
as follows~: the pairing $A \{\ell \} \times B^{(\ell)} \to \Q/\Z$ 
has trivial right kernel and divisible left kernel.
}
\end{rema}

The following theorem extends the function field case of 
\cite{MilADT}, Th. II.4.6. (which corresponds to $\C=\T$ or 
$\C=\T[1]$, where $\T$ is a torus).

\begin{theo} \label{refin2}

\smallskip

a) The pairing (\ref{leftpairing})
induces a perfect duality between the discrete torsion group
$H^3_c(U,\widehat {\mathcal C})$ and the finite-type $\widehat \Z$-module
$H^{-1}(U,{\mathcal C})_{\wedge}$, resp. between 
the discrete torsion group
$H^2_c(U,\widehat {\mathcal C})$ and the finite-type $\widehat \Z$-module
$H^0(U,{\mathcal C})_{\wedge}$.

\smallskip

b) Assume $U \neq X$. The pairing (\ref{leftpairing})
induces a perfect duality between the discrete torsion group
$H^1(U,{\mathcal C})$
and the finite-type $\widehat \Z$-module
$H^1_c(U,\widehat {\mathcal C})_{\wedge}$, resp. between
the discrete torsion group
$H^2(U,{\mathcal C})$ and the finite-type $\widehat \Z$-module
$H^0_c(U,\widehat {\mathcal C})_{\wedge}$.

\end{theo}

\dem{} a) Let $\ell$ be a prime number. Then the map $\psi'$ 
of Lemma~\ref{avlem} a) is an isomorphism for $i=-2$, 
which yields the first point (recall that $H^3_c(U,\widehat {\mathcal C})$
and $H^2_c(U,\widehat {\mathcal C})$ are torsion of cofinite type
by Proposition~\ref{finivanish1}; also $H^{-1}(U,\C)$ and $H^0(U,\C)$ 
are finitely generated by Proposition~\ref{finityp}).

\smallskip

In the case $U=X$, 
the second point is a duality between $H^2(X,\widehat {\mathcal C})$
and $H^0(X,C)_{\wedge}$, which follows from Lemma~\ref{avlem} b)
in the case $i=1$. 
Now assume $U \neq X$ and let ${\mathcal T}$ be a $U$-torus.
By the first point applied
to $\mathcal C=\mathcal T [1]$,
the group $H^3_c(U,\widehat {\mathcal T})$ is dual to
$H^0(U,\mathcal T)_{\wedge}$. By Lemma~\ref{avlem} a) with $i=-1$
and $\mathcal C=\mathcal T [1]$, the finite group $H^2_c(U,\widehat \T)$ 
(cf. Lemma~\ref{h1lem}, c)
is dual to the finite group $H^1(U,T)$ (cf. Lemma~\ref{fini1}, a): indeed
a finite group coincides
with its $\ell$-adic completion and doesn't contain a non trivial
divisible subgroup.
There is a commutative diagram with exact lines (observe that the 
second line is obtained by applying the $\ell$-completion functor
to an exact sequence of finitely generated group, then
dualizing an exact sequence of profinite groups):

{\small
$$
\begin{CD}
H^2_c(U,\widehat \T_2)\{ \ell \} @>>> H^2_c(U,\widehat
\T_2)\{ \ell \}
@>>> H^2_c(U,\widehat \C)\{\ell \} @>>> H^3_c(U,\widehat
\T_2)\{ \ell \}
@>>> H^3_c(U,\widehat \T_1)\{ \ell \} \cr
@VV{f_1}V @VV{f_2}V @VV{h}V @VV{g_1}V @VV{g_2}V \cr
(H^1(U,{\mathcal T_2})^{(\ell)})^* @>>>  (H^1(U,
{\mathcal T_1})^{(\ell)})^* @>>> (H^0(U,{\mathcal C})^{(\ell)})^*
@>>> (H^0(U,{\mathcal T_2})^{(\ell)})^* @>>> (H^0(U,{\mathcal T_1})^{(\ell)})^*.
\end{CD}
$$
}

Now $h$ is an isomorphism by the five-lemma, whence the result.

\smallskip

b) Consider diagram (\ref{finidiag1}) for $i=0$. 
By a), the left vertical map is an isomorphism
and by Proposition~\ref{avfini} 1., 
the middle vertical map is an isomorphism, hence
$\psi_m$ is an isomorphism from $_{\ell^m} H^1(U, {\mathcal C})$
to $(H^1_c(U,\widehat \C)/\ell^m)^*$. Taking direct limit over $m$, then
direct sum over all prime $\ell$, yields the duality between the torsion 
group of cofinite type (cf. Proposition~\ref{finivanish1})
$H^1(U, {\mathcal C})$ and the finite type $\widehat \Z$-module 
(cf. Proposition~\ref{finityp}) $H^1_c(U,\widehat \C)_{\wedge}$. 

\smallskip

Lemma~\ref{avlem} a) for $i=1$ yields
that for $U \neq X$, the map 
$\psi$ an isomorphism, which immediately gives the duality between 
the torsion group of finite type (cf. Proposition~\ref{finivanish1}) 
$H^2(U,\C)$ and the finite type $\widehat \Z$-module (cf. 
Proposition~\ref{finityp}) $H^0_c(U, \widehat \C)_{\wedge}$. 

\enddem 

\begin{rema}
{\rm In the case $U=X$, the first assertion of Theorem~\ref{refin2} 
b) should be replaced by a duality between 
$H^1(X,\widehat \C)$ and $H^1(X, {\mathcal C})_{\wedge}$ (see 
Theorem~\ref{refin1} b) below in the case $U=X$). The second assertion  
(duality 
between $H^2(X,\C)$ and $H^0(X, \widehat \C)_{\wedge}$) actually still
holds, cf. Theorem~\ref{refin1} a).

}
\end{rema}

The following duality theorem has the same flavour as \cite{MilADT}, 
Th~II.4.6. b) (but one should be careful that in the number field case, 
the case $r=3$ of the latter does not hold in general, see also
Remark~\ref{cdn2}).

\begin{theo} \label{refin1}

\smallskip

a) The pairing (\ref{rightpairing})
induces a perfect duality between the discrete torsion group
$H^3_c(U,\mathcal C)$ and the finite-type $\widehat \Z$-module
$H^{-1}(U,\widehat {\mathcal C})_{\wedge}$, resp. between
the discrete torsion group
$H^2_c(U,\mathcal C)$ and the finite-type $\widehat \Z$-module
$H^0(U,\widehat {\mathcal C})_{\wedge}$.

\smallskip

b) The pairing (\ref{rightpairing})
induces a perfect duality between the discrete torsion group
$H^1(U,\widehat {\mathcal C})$ and the profinite
group $H^1_c(U,\mathcal C)_{\wedge}$, resp. between
the discrete torsion group
$H^2(U,\widehat {\mathcal C})$ and the profinite
group $H^0_c(U,\mathcal C)_{\wedge}$.

\end{theo}

\dem{} a) Let $\ell$ be any prime number. The map $\varphi'$ of
Lemma~\ref{avlem} b) is an isomorphism for $i=-2$, which yields
the first point (Proposition~\ref{avfini} 2. yields that 
$H^3_c(U,\C)$ and $H^2_c(U,\C)$ are torsion groups of cofinite type; 
Proposition~\ref{finityp} gives that $H^{-1}(U,\widehat \C)$ and 
$H^0(U, \widehat \C)$ are finitely generated).

\smallskip

There is a commutative diagram with exact lines
{\small

$$
\begin{CD}
H^2_c(U,\mathcal T_1)\{ \ell \} @>>> H^2_c(U,\mathcal T_2)\{ \ell \}
@>>> H^2_c(U,\mathcal C)\{\ell \} @>>> H^3_c(U,\mathcal T_1)\{ \ell \}
@>>> H^3_c(U,\mathcal T_2)\{ \ell \} \cr
@VV{f_1}V @VV{f_2}V @VV{h}V @VV{g_1}V @VV{g_2}V \cr
(H^1(U,\widehat {\mathcal T_1})^{(\ell)})^* @>>>  (H^1(U,\widehat
{\mathcal T_2})^{(\ell)})^* @>>> (H^0(U,\widehat {\mathcal C})^{(\ell)})^*
@>>> (H^0(U,\widehat {\mathcal T_1})^{(\ell)})^* @>>> (H^0(U,\widehat {\mathcal T_2})^{(\ell)})^*.
\end{CD}
$$
}

The maps $g_1$ and $g_2$ are isomorphisms by the first point
applied to
$\mathcal C=\mathcal T_1$, $\mathcal C=\mathcal T_2$. The maps
$f_1$ and $f_2$ are isomorphisms by Lemma~\ref{avlem} b) applied 
to the same complexes (map $\varphi$ in the case $i=-1$): indeed
for a $U$-torus $\mathcal T$, the groups
$H^1(U,\widehat {\mathcal T})$
and $H^2_c(U,{\mathcal T})$ are finite
(Lemma~\ref{h1lem} a) and b), hence they coincide with their
$\ell$-adic completions and do not contain a non trivial divisible 
subgroup.
Therefore $h$ is an isomorphism by the five-lemma, whence the second 
point.

\smallskip

b) Consider diagram (\ref{finidiag3}) for $i=0$. By
a), the left vertical map is an isomorphism
and the middle vertical map is an isomorphism by
Proposition~\ref{avfini} 1., hence
$\varphi_m$ is an isomorphism from $_{\ell^m} H^1(U, \widehat {\mathcal C})$
to $(H^1_c(U,\mathcal C)/\ell^m)^*$. Taking direct limit over $m$, then
direct sum over all prime $\ell$, yields the duality between
$H^1(U, \widehat {\mathcal C})$ (which is torsion by 
Proposition~\ref{finivanish2}, but not necessarily of cofinite
type, cf. Remark~\ref{nocofini}) and $H^1_c(U,\mathcal C)_{\wedge}$.

\smallskip

Now consider diagram (\ref{finidiag3}) for $i=1$. By the previous duality, 
the left vertical map induces an isomorphism between 
$H^1(U, \widehat \C)/\ell^m$ and $(_{\ell^m} (H^1_c(U,\C)_{\wedge}))^*$. 
Since $H^1_c(U,\C)$ is in the class $\mathcal E$ (that is:
it is the product of a finite type group by a profinite 
group) by Proposition~\ref{injectcompl}, 
the $\ell^m$-torsion of $H^1_c(U,\C)$ and of 
$H^1_c(U,\C)_{\wedge}$ coincide, hence the left vertical map is 
actually an isomorphism and the right vertical map $\psi_m$ is an 
isomorphism as well (the middle vertical map
is an isomorphism by Proposition~\ref{avfini} 2.).
Taking direct limit and direct sum over 
all prime $\ell$, we get 
the duality between the torsion group (cf. Proposition~\ref{finivanish2}) 
$H^2(U,\widehat \C)$ and 
the profinite group $H^0_c(U,\C)_{\wedge}$.

\enddem

\begin{prop} \label{duald1}
The pairing (\ref{leftpairing}) for $r=1$ 
induce a perfect pairing of finite groups 
\begin{equation} \label{d1pairing}
D^1(U, \mathcal C) \times D^1(U,\widehat {\mathcal C}) \to \Q/\Z.
\end{equation}
\end{prop}

\dem{}  Fix a prime number $\ell$. 
There is a commutative diagram~: 
\begin{equation} \label{diaglocglob}
\begin{CD}
0 @>>> D^1(U,\mathcal C)\{ \ell \} @>>> H^1(U,\mathcal C)\{ \ell \} @>>> 
\bigoplus_{v \not \in U} H^1(K_v,C)\{ \ell \} \cr
&& && @VV{\psi}V @VV{\beta}V \cr
0 @>>>  (D^1(U,\widehat {\mathcal C}) ^{(\ell)})^* @>>>
 (H^1_c(U, \widehat {\mathcal C}) ^{(\ell)})^* @>>> (\bigoplus_{v \not \in U} 
H^0(K_v,\widehat C)^{(\ell)})^* .
\end{CD}
\end{equation}
The first line is exact by definition of $D^1(U,\mathcal C)$. 
The sequence 
$$\bigoplus_{v \not \in U} H^0(K_v,\widehat C) \to 
H^1_c(U, \widehat {\mathcal C}) \to D^1(U,\widehat {\mathcal C}) \to 0$$
is also exact by definition of $D^1(U,\widehat {\mathcal C})$.
Using Lemma~\ref{exactcomp} and the fact that the
$\ell$-adic completion functor $^{(\ell)}$ commutes 
with finite direct sums, the sequence 
$$\bigoplus_{v \not \in U} H^0(K_v,\widehat C)^{(\ell)} \to 
H^1_c(U, \widehat {\mathcal C})^{(\ell)} \to 
D^1(U,\widehat {\mathcal C})^{(\ell)} \to 0$$
of profinite groups is exact as well, and its dual sequence (which is the
second line of the diagram) remains exact.

\smallskip

The commutative diagram (\ref{diaglocglob}) 
defines a map $$\theta : D^1(U,\mathcal C)\{ \ell \} \to 
(D^1(U,\widehat {\mathcal C}) ^{(\ell)})^*.$$ We observe that 
by \cite{demphd} Th.~3.1 (which is the local duality theorem), 
the map $\beta$ is an isomorphism, 
and we also know by Theorem~\ref{refin2}  
that $\psi$ is an isomorphism.
By diagram chasing $\theta$ is an isomorphism. 
Since this holds for every
prime $\ell$ (including $\ell=p$), the proposition is proven, the finiteness
of $D^1(U,\mathcal C)$ and $D^1(U,\widehat {\mathcal C})$ being known by  
Proposition~\ref{finityp} c).

\enddem

\begin{rema}
{\rm Of course, the pairing (\ref{d1pairing}) can also be defined 
via the pairing (\ref{rightpairing}).
}
\end{rema}

\begin{lem} \label{d2lem}
Assume $U \neq X$. Then
$D^2(U,\mathcal C)$ and $D^0(U, \widehat \C)$ are finite.
\end{lem}

\dem{} Using the exact triangle (\ref{devistm}) and the fact that
$\cok \rho:=\mathcal T$ is a torus, we know that
$D^2(U,\mathcal T)$ is finite and is sufficient to show that
$H^3(U,\ker \rho)$ is finite to get the finiteness of $D^2(U,\mathcal C)$.
But $\ker \rho$ is a group of multiplicative
type, so there is an exact sequence
$$0 \to \mathcal T_1 \to \ker \rho \to \mathcal F \to 0, $$
where $\mathcal F$ is a finite group of multiplicative type and
$\mathcal T_1$ is a torus.
Since $H^3(U,\mathcal T_1)=0$ by Proposition~\ref{finivanish1} a) and
$H^3(U,\mathcal F)=0$ (cf. Remark~\ref{fini2};
it is dual to $H^0_c(U,\widehat {\mathcal F})$,
which is zero because $U \neq X$), the group $H^3(U,\ker \rho)$
is actually zero.

\smallskip

The group $D^0(U,\widehat {\ker \rho})$ is trivial thanks to the 
assumption $U \neq X$. Thus the exact triangle (\ref{devistmbis}) 
shows that $D^0(U,\widehat \C)$ is finite because 
so is $H^1(U,\widehat \T)$ (Lemma~\ref{h1lem} a).

\enddem

\begin{lem} \label{d0lem}
Assume $U \neq X$. Then the groups $D^0(U,\C)$ and $D^2(U,\widehat \C)$ 
are finite. 
\end{lem}

\dem{} For a $U$-torus $\T$, we have $D^0(U,\T)=0$ because $U \neq X$. 
For a $U$-group of multiplicative type $\mathcal M$, we also know 
(Remark~\ref{fini2}) that $H^1(U,\mathcal M)$ is finite, whence 
the finiteness of $D^0(U,\C)$ via the 
exact triangle (\ref{devistm}). 

\smallskip 

Exact triangle (\ref{devistmbis}) and the vanishing of $H^3(U, \widehat \T)$ 
for a $U$-torus $\T$ (Lemma~\ref{h1lem} d) imply that $D^2(U,\widehat \T)$
injects into $D^2(U, \widehat {\mathcal M})$, so it only remains to show
that the latter is finite. We show that $H^2_c(U,\widehat 
{\mathcal M})$ is finite. By d\'evissage it is sufficient to prove this 
when $M$ is a finite group of multiplicative type and when $M$ is 
a torus. The first case follows from Lemma~\ref{fini1} b) and the second one 
from Lemma~\ref{h1lem} c).

\enddem

\begin{prop} \label{duald2}
The pairing (\ref{leftpairing}) for $r=2$
induce a perfect pairing of finite groups
\begin{equation} \label{d2pairing}
D^2(U, \mathcal C) \times D^0(U,\widehat {\mathcal C}) \to \Q/\Z.
\end{equation}
\end{prop}

\dem{} The argument is exactly the same as in the proof of 
Proposition~\ref{duald1}, using now the commutative diagram with 
exact lines: 
\begin{equation} \label{diaglocglob2}
\begin{CD}
0 @>>> D^2(U,\mathcal C)\{ \ell \} @>>> H^2(U,\mathcal C)\{ \ell \} @>>>
\bigoplus_{v \not \in U} H^2(K_v,C)\{ \ell \} \cr
&& && @VVV @VVV \cr
0 @>>>  (D^0(U,\widehat {\mathcal C}) ^{(\ell)})^* @>>>
 (H^0_c(U, \widehat {\mathcal C}) ^{(\ell)})^* @>>> (\bigoplus_{v \not \in U}
H^{-1}(K_v,\widehat C)^{(\ell)})^* .
\end{CD}
\end{equation}

Indeed the right vertical map is an isomorphism by \cite{demphd}, Th.~3.1
and the middle vertical map is an isomorphism as well by 
Theorem~\ref{refin2}, b). It remains to apply Lemma~\ref{d2lem}.

\enddem

\begin{prop} \label{duald0}
The pairing (\ref{rightpairing}) for $r=0$
induce a perfect pairing of finite groups
\begin{equation} \label{d0pairing}
D^0(U, \mathcal C) \times D^2(U,\widehat {\mathcal C}) \to \Q/\Z.
\end{equation}
\end{prop}

\dem{} Again the argument is the same as in 
Proposition~\ref{duald1}, using this time the commutative diagram with
exact lines:
\begin{equation} \label{diaglocglob0}
\begin{CD}
0 @>>> D^2(U,\widehat C)\{ \ell \} @>>> H^2(U,\widehat C)\{ \ell \} @>>>
\bigoplus_{v \not \in U} H^2(K_v,\widehat C)\{ \ell \} \cr
&& && @VVV @VVV \cr
0 @>>>  (D^0(U,{\mathcal C}) ^{(\ell)})^* @>>>
 (H^0_c(U, {\mathcal C}) ^{(\ell)})^* @>>> (\bigoplus_{v \not \in U}
H^{-1}(K_v,C)^{(\ell)})^* .
\end{CD}
\end{equation}

The right vertical map is an isomorphism by \cite{demphd}, Th.~3.1
(recall that by Proposition~\ref{injectcompl}, the groups
$H^{-1}(K_v,C)_{\wedge}$ and $H^{-1}(K_v,C)^{\wedge}$ coincide),
and the middle vertical map is an isomorphism as well by
Theorem~\ref{refin1}, b); Lemma~\ref{d0lem} then yields the result.

\enddem

\begin{rema} \label{cdn2}
{\rm Again there are analogous results over a number field~: 

\begin{itemize}

\item Proposition~\ref{avfini} 1. holds except that for $n$ even, 
the vanishing statements for $i \geq 2$ do not hold anymore
if $\Omega_{\RR} \neq \emptyset$. In Proposition~\ref{avfini} 2.,
all groups involved are finite, but the vanishing statements 
for $i \leq -2$ are in general false if $n$ is even and
$\Omega_{\RR} \neq \emptyset$.

\item In Proposition~\ref{finivanish1} and Proposition~\ref{finivanish2}, 
the vanishing statements must be replaced by finiteness statements 
if $\Omega_{\RR} \neq \emptyset$  for the following groups: 
$H^i(U,\C)$ for $i \geq 4$, $H^i_c(U,\widehat \C)$ for $i \leq -2$,
$H^i(U,\widehat \C)$ for $i \geq 4$, $H^i_c(U,\C)$ for $i \leq -2$. 
The groups $H^3(U,\C)$, $H^{-1}_c(U,\widehat \C)$, and 
$H^{-1}_c(U,C)$ are still finite if $U \neq X$ (resp. zero if 
$U \neq X$ and $\Omega_{\RR} = \emptyset$).
Also, the group $H^1(U,\C)$ is finite even if $U=X$ and 
the finiteness assertions in Proposition~\ref{finivanish2} a) 
hold without any condition on $i$.
Finally, 
{\it the vanishing of $H^3(U,\widehat \C)$ does not hold anymore}
(see Remark~\ref{cdn1} about $H^3(U, \widehat \T)$) even 
for $\Omega_{\RR}=\emptyset$.

\item Lemma~\ref{avlem} is unchanged, except that the restriction 
$U \neq X$ can be removed in a) for $i=1$. 

\item Theorem~\ref{refin2} is unchanged (which gives a more precise statement 
than \cite{demphd}, Th.~4.3) except that the assumption 
$U \neq X$ can be removed in b). Theorem~\ref{refin1} a) is still true,
as is the first assertion of Theorem~\ref{refin1} b), but 
{\it not the second assertion of Theorem~\ref{refin1} b)}: the 
pairing $H^2(U,\widehat \C) \{ \ell  \} \times H^0_c(U,\C)^{(\ell)}$ 
has trivial right kernel and divisible left kernel, but 
for triviality of the left kernel we need $\ell$ invertible 
on $U$ and Leopoldt's conjecture. 

\item Proposition~\ref{duald1} is unchanged (this removes the condition $\ell 
\in \calo_U^*$ in \cite{demphd}, Cor~4.7). 
Lemma~\ref{d2lem} also holds
(in the proof, the groups $H^3(U,\mathcal T_1)$ and 
$H^0_c(U,\widehat {\mathcal F})$ might be only finite if
$\Omega_{\RR} \neq \emptyset$, but this does not affect the result),
as does Proposition~\ref{duald2} (the assumption $\ker \rho$ finite
made in \cite{demphd}, Lemma~5.13, is not necessary).
The first part of Lemma~\ref{d0lem} still holds,
but {\it not its second part} because in general the $\ell$-primary part
of $H^3(U,\widehat \T)$ is infinite if $\ell$ is not invertible 
on $U$ (and even for $\ell \in \calo_U^*$, the finiteness of 
$H^3(U,\widehat \T)\{ \ell \}$ relies on Leopoldt's conjecture).
Similarly  {\it Proposition~\ref{duald0} does not hold anymore 
in general}, we only get that the pairing (\ref{d0pairing}) 
has trivial left kernel and 
divisible right kernel (see also \cite{demphd}, paragraph~5.4. for 
a variant).

\end{itemize}

}
\end{rema}

\section{Poitou-Tate exact sequences} \label{three}

Let $C=[T_1 \to T_2]$ be a complex of $K$-tori with dual 
$\widehat C=[\widehat T_2 \to \widehat T_1]$. We can choose a non empty 
Zariski open subset $U_0$ of $X$ such that $C$ extends to a complex 
$\mathcal C=[{\mathcal T}_1 \stackrel{\rho}{\to}{\mathcal T}_2]$
of $U_0$-tori with dual 
$\widehat {\mathcal C}$. For every integer $i$ and every $K$-group 
scheme (or bounded complex of $K$-group schemes) 
$M$ (e.g. $M=T$, $M=\widehat T$), define
$$\Sha^i(M):={\rm Ker} \, [H^i(K,M) \to \prod_{v \in X^{(1)}} H^i(K_v,M)].$$

\begin{lem} \label{shad1}

There exists a non empty Zariski open subset $U_1 \subset U_0$ 
such that for every Zariski open subset $V \subset U_1$:

\smallskip

a)  For $i \in \{1, 2 \}$, 
the restriction map $r_{U_1,V} : H^i(U_1,\mathcal C) \to H^i(V,\mathcal C)$ 
induces isomorphisms 
$$D^i(U_1,\mathcal C) \cong D^i(V,\mathcal C) \cong \Sha^i(C).$$

\smallskip

b) For $r \in \{ 0, 1\}$, 
the canonical map $H^r(V, \widehat {\mathcal C}) \to H^r(K,
{\widehat \C})$ is injective and identifies $D^r(V,\widehat {\mathcal C})$
with $\Sha^r(\widehat C)$.

\end{lem}

\dem{} We can deal with the two properties a) and b) separately 
(up to taking the intersections of the various provided $U_1$).

\smallskip

a) Let's start with arbitrary non empty Zariski open subsets 
$V \subset U \subset U_0$. Take $i \in \{1 ,2 \}$. 
For all $v \in U$, we have $H^i(\calo_v,\mathcal C)=0$ by 
Lemma~\ref{vanishlocal}, which implies that 
the image of $D^i(U,\mathcal C)$ by $r_{U,V}$ is contained in  
$D^i(V,\mathcal C)$. The induced map $D^i(U,\mathcal C) \to 
D^i(V,\mathcal C)$ is surjective thanks to the compatibility of 
the covariant map $H^i_c(V,\mathcal C) \to H^i_c(U,\mathcal C)$ 
with $r_{U,V}$ (\cite{DemH} Prop~2.1 (3)). Since all $D^i(U,\mathcal C)$ 
are finite by Proposition~\ref{finityp} c) and Lemma~\ref{d2lem},
the decreasing sequence 
of positive integers $\# D^i(U,\mathcal C)$ (when $U$ becomes 
smaller and smaller) must stabilize for some $U=U_1$. We get 
an isomorphism from $D^i(U_1,\mathcal C)$ to $D^i(V,\mathcal C)$
for all $V \subset U_1$. Since 
$H^i(K,C)$ is the direct limit over $V$ of the $H^i(V,\mathcal C)$, we 
get an injective map $u: D^i(U_1,\mathcal C) \to H^i(K,C)$. As 
$D^i(U_1,\mathcal C)$ is the same as $D^i(V,\mathcal C)$ for every 
$V \subset U_1$, the image of $u$ is contained in $\Sha^i(C)$ (because 
its restriction to $H^i(K_v,C)$ is zero for all $v \not \in V$ and 
$V$ can be taken arbitrarily small). Conversely, every element 
of $\Sha^i(C)$ can be lifted to an $a \in H^i(V, \mathcal C)$ for some 
$V$, and by definition $a \in D^i(V,\mathcal C)=D^i(U_1,\mathcal C)$, so 
the image of $u$ contains $\Sha^i(C)$.

\smallskip

b) Let $V \subset U_0$ be an arbitrary non empty Zariski 
open subset. Let $r \in \{0, 1 \}$.
The injectivity 
of $H^r(V, \widehat {\mathcal C}) \to H^r(K,\widehat C)$ has been 
proven in Proposition~\ref{finityp} a).
Identifying 
now $D^r(V,\widehat {\mathcal C})$ with a subgroup of 
$H^r(K,\widehat C)$, we get (again using the maps
$H^r(V,\widehat {\mathcal C}) \to H^r(U,\widehat {\mathcal C})$ for 
$V \subset U \subset U_0$) a decreasing sequence of 
finite subgroups (when $V$ becomes smaller and smaller), which 
stabilizes for some $U_1$. Since $D^r(U_1,\widehat {\mathcal C})$ is 
also $D^r(V,\widehat {\mathcal C})$ for every $V \subset U_1$, we
have $D^r(U_1,\widehat {\mathcal C}) \subset \Sha^r(\widehat C)$. 
On the other hand, every element of $\Sha^r(\widehat C)$ comes from
$H^r(V,\widehat {\mathcal C})$ for some $V \subset U_1$, and it is then 
automatically in $D^r(V,\widehat {\mathcal C})=D^r(U_1,\widehat {\mathcal C})$
because it is everywhere locally trivial. 

\enddem

\begin{theo} \label{sha12theo}
There are perfect pairing of finite groups 
$$\Sha^1(C) \times \Sha^1(\widehat C) \to \Q/\Z .$$
$$\Sha^2(C) \times \Sha^0(\widehat C) \to \Q/\Z .$$
\end{theo}

\dem{} This follows immediately from Lemma~\ref{shad1} and 
Proposition~\ref{duald1} (resp. Proposition~\ref{duald2}) applied to $U_1$.

\enddem

\begin{lem} \label{injecth0}
There exists a non empty Zariski open subset $U_1$ of $U_0$ such that 
for every non empty Zariski open subset $V$ of $U_1$:

\smallskip

-the restriction map
$H^0(V,\C) \to H^0(K,C)$ is injective.

\smallskip

-For all non empty Zariski open subsets $W \subset V$, the canonical map 
$$j_{W,V} : H^2_c(W,\widehat \C) \to H^2_c(V,\widehat \C)$$ is surjective and 
the image of $D^2(V,\widehat \C)$ by the restriction map  
$$r_{V,W} : H^2(V,\widehat \C) \to H^2(W,\widehat \C)$$ is a subgroup of 
$D^2(W,\widehat \C)$.
\end{lem}

\dem{} Let $U \subset U_0$ be a non empty Zariski open subset. 
By the exact triangle (\ref{devistm}), there is a commutative 
diagram with exact lines

$$ 
\begin{CD}
0 @>>> H^1(U,\mathcal M) @>>> H^0(U,\C) @>>> H^0(U,\T_0) \cr
&& @VVV @VVV @VVV \cr
0 @>>> H^1(K,M) @>>> H^0(K,C) @>>> H^0(K,T_0) ,
\end{CD}
$$
where $\mathcal M$ is a $U$-group of multiplicative type with generic 
fibre $M$ and $\T_0$ is a $U$-torus.
Since the right vertical map is clearly injective, it is sufficient to 
prove the injectivity of the left vertical map for $U$ small enough.
We can write $\mathcal M$ as an extension 
$$0 \to \T \to \mathcal M \to \mathcal F \to 0$$
of a finite $U$-group of multiplicative type $\mathcal F$ by a 
$U$-torus $\T$. This yields a commutative diagram with exact lines
$$ 
\begin{CD}
0 @>>> H^0(U,\mathcal F) @>>> H^1(U,\T) @>>> H^1(U,\mathcal M) @>>> 
H^1(U,\mathcal F) \cr
&& @VVV @VVV @VVV @VVV \cr
0 @>>> H^0(K,F) @>>> H^1(K,T) @>>> H^1(K,M) @>>> 
H^1(K,F).
\end{CD}
$$
Since $\mathcal F$ is finite (hence proper) over $U$, the left vertical 
map is an isomorphism and the right vertical map is injective. 
It is therefore sufficient to prove that for a $U$-torus $\T$, the 
restriction map $H^1(U,\T) \to H^1(K,T)$ is injective for $U$ sufficiently 
small. 
Set $N_U=\ker [H^1(U,\T) \to H^1(K,T)]$. For every Zariski open 
subset $V \to U$, the restriction map 
$H^1(U,\T) \to H^1(V,\T)$ induces a homomorphism $i_{U,V} : N_U \to N_V$.
By Lemma~\ref{uvlem},
this homomorphism is surjective. Lemma~\ref{fini1} implies that 
the group $N_V$ is finite, and the decreasing sequence of positive 
integers $(\# N_V), V \subset U_0$ must stabilize for some
$V=U_1 \subset U_0$. Then the maps $i_{U_1,V}$ for $V \subset U_1$
are isomorphisms,
which implies (passing to the limit) that the restriction map 
$N_{U_1} \to H^1(K,T)$ is injective. By definition of $N_{U_1}$, this 
means that $N_{U_1}=0$, hence $N_V=0$ for every $V \subset U_1$. This gives
the first point. 

\smallskip 

For $W \subset V \subset U_1$, 
the restriction map $H^0(V,\C) \to H^0(W,\C)$ is injective
because so is its composition with $H^0(W,\C) \to H^0(K,\C)$. As 
$H^0(V,\C)$ and  $H^0(W,\C)$ are finitely generated by
Proposition~\ref{finityp} b), the induced map 
$H^0(V,\C)_{\wedge} \to H^0(W,\C)_{\wedge}$ is still injective. 
By Theorem~\ref{refin2}, the dual map $H^2_c(W,\widehat \C) \to 
H^2_c(V,\widehat \C)$ is surjective. Now the compatibility 
of $r_{V,W}$ with $j_{W,V}$ (\cite{DemH}, Prop~2.1 c) gives that
$r_{V,W} (D^2(V,\widehat \C)) \subset D^2(W,\widehat \C)$.

\enddem

\begin{theo} \label{sha0theo}
There is a perfect pairing of finite groups
$$\Sha^0(C) \times \Sha^2(\widehat C) \to \Q/\Z .$$
\end{theo}

\dem{} As in the proof of Lemma~\ref{shad1}, Lemma~\ref{injecth0} and 
Lemma~\ref{d0lem} imply that for a sufficiently small Zariski open subset 
$U \subset U_0$, we have $\Sha^0(\C)=D^0(U,\C)$ and $\Sha^2(\widehat \C) \cong
D^2(U,\widehat \C)$. Now apply Proposition~\ref{duald0}.

\enddem

For each integer $i$, denote by $\prod^{'}_{v \in X^{(1)}} 
H^i(K_v,C)$ (resp. $\prod^{'}_{v \in X^{(1)}} 
H^i(K_v,C)_{\wedge}$)
the restricted product of the $H^i(K_v,C)$ (resp. of the $H^i(K_v,C)_{\wedge}$)
with respect 
to the $H^i_{\nr}(K_v,C)$ (resp. to the image of $H^1(\calo_v,\C)$ in 
$H^i(K_v,C)_{\wedge}$). The same notation stands for $\widehat C$.
The groups $\prod^{'}_{v \in X^{(1)}} H^i(K_v,C)$ and 
$\prod^{'}_{v \in X^{(1)}} H^i(K_v,\widehat C)$ are equipped with their 
restricted product topology (associated to the topology previously 
defined on the $H^i(K_v,C)$ and $H^i(K_v,\widehat C)$). All groups 
$H^i(K,C)$ (resp. $H^i(K,\widehat C)$) are equipped with the discrete 
topology. 

\begin{lem} \label{discretlem}
Let $i$ be an integer. Then the image of $H^i(K,C)$
in $\prod^{'}_{v \in X^{(1)}} H^i(K_v,C)$ is discrete for the subspace 
topology. The same holds if $C$ is replaced by $\widehat C$.
\end{lem}

\dem{} As the local fields $K_v$ are of strict cohomological dimension $2$, 
the statement is obvious except for $-1 \leq i \leq 2$.
Fix a Zariski open subset $U \subset U_0$ with $U \neq X$.
All groups $H^i(K_v,\widehat C)$ are discrete, so the subgroup 
$E:=\prod_{v \not \in U} \{ 0 \} \times \prod_{v \in U}
H^i_{\nr}(K_v, \widehat C)$ is open in $\prod^{'}_{v \in X^{(1)}}
H^i(K_v,\widehat C)$. Let $I$ be the image of
$H^i(K,\widehat C)$ in $\prod^{'}_{v \in X^{(1)}} H^i(K_v,\widehat C)$.
Every element of $H^1(K,\widehat C)$ comes from 
$H^1(V, \widehat \C)$ for some $V \subset U$, hence by 
Lemma~\ref{uvlem}, there is a surjection $D^i(U,\widehat \C) \to I \cap E$.
Since all groups $D^i(U,\widehat \C)$ are finite by
Proposition~\ref{finityp} c), Lemma~\ref{d2lem} and Lemma~\ref{d0lem},
this implies that $I \cap E$ is finite, hence $I$ is discrete. 

\smallskip

The same argument shows that the image $J$ of $H^i(K,C)$
in $\prod^{'}_{v \in X^{(1)}} H^i(K_v,C)$ is discrete for $i \geq 1$.
For $i \in \{-1, 0 \}$,
this is an immediate consequence of Lemma~\ref{closedimage} (again combined 
with Lemma~\ref{uvlem}). 

\enddem

\begin{lem} \label{ptstart}
Let $U \subset U_0$ be a non empty Zariski open subset with $U \neq X$.

\smallskip

a) There are exact sequences 
$$H^0(U,\C) \to \prod_{v \not \in U} H^0(K_v,C) \times \prod_{v \in U} 
H^0_{\nr}(K_v,C) \to H^1(K,\widehat C)^* .$$
$$0 \to 
H^{-1}(U,\C) \to \prod_{v \not \in U} H^{-1}(K_v,C) \times \prod_{v \in U}
H^{-1}_{\nr}(K_v,C) \to H^2(K,\widehat C)^* .$$

\smallskip

b) There are exact sequences 
$$H^2(U,\widehat \C) \to \prod_{v \not \in U} H^2(K_v,\widehat C)
\times \prod_{v \in U} 
H^2_{\nr}(K_v,\widehat C) \to H^{-1}(K,C)^* \to 0 .$$
$$ 
H^{1}(U,\widehat \C) \to \prod_{v \not \in U} H^{1}(K_v,\widehat C)
\times \prod_{v \in U} H^{1}_{\nr}(K_v,\widehat C) \to H^0(K,C)^* \to
D^2(U,\widehat C).$$

\end{lem}

\dem{} a) Let $V \subset U$ be a non empty Zariski open subset. Let
$i \in \{-1, 0 \}$. By Lemma~\ref{uvlem}, 
we have an exact sequence 
$$ H^i(U,\C) \to \prod_{v \not \in U} H^i(K_v,C) \times
\prod_{v \in U-V} H^i_{\nr}(K_v,C) \to H^{i+1}_c(V,\C).$$

By Proposition~\ref{injectcompl}, the map $H^{i+1}_c(V,\C) \to
H^{i+1}_c(V,\C)_{\wedge}$ 
is injective, thus by Theorem~\ref{refin1} we get an exact sequence 
$$H^i(U,\C) \to \prod_{v \not \in U} H^i(K_v,C) \times 
\prod_{v \in U-V} H^i_{\nr}(K_v,C) \to H^{1-i}(V,\widehat \C)^*,$$
where $H^{1-i}(V,\widehat \C)$ is a discrete torsion group.
Besides, the kernel of the first map 
is a subgroup of $D^i(U,\C)$, hence it is finite for $i=0$
by Lemma~\ref{d0lem}. This kernel is also obviously zero for $i=-1$ as soon 
as $V \neq U$. 
This implies that the inverse limit of this exact sequence 
(when $V$ runs over all non empty Zariski open subsets of $U$) remains exact,
which yields the result.

\smallskip


\smallskip

b) We apply again Lemma~\ref{uvlem} and observe that for $i \in \{1, 2 \}$:

\smallskip

-we have $H^{i+1}_c(V,\widehat \C) \simeq 
(H^{1-i}(V,\C)_{\wedge})^* \simeq H^{1-i}(V,\C)^*$ by Theorem~\ref{refin2}, 
because the discrete finitely generated (cf. Proposition~\ref{finityp}, b) 
group $H^{1-i}(V,\C)$ and its completion have same dual.

-the groups $D^i(U,\widehat \C)$ are finite (Lemma~\ref{d0lem} and 
Proposition~\ref{finityp} c).

\smallskip

Now the same method as in a) gives the exactness of 
$$H^i(U,\widehat \C) \to 
\prod_{v \not \in U} H^i(K_v,\widehat C) \times \prod_{v \in U}
H^i_{\nr}(K_v,\widehat C) \to H^{1-i}(K,C)^*.$$ 

\smallskip

Besides, by \cite{DemH}, Prop~2.1, there is a commutative diagram 
with exact lines:

$$\begin{CD} 
&&&& D^{i+1}(U,\widehat \C)  \cr 
&&&& @AAA \cr
\prod_{v \in U-V} H^i_{\nr}(K_v,\widehat C) @>>> 
H^{i+1}_c(V,\widehat \C) @>>> H^{i+1}_c(U,\widehat \C) \cr
&& @AAA @AAA \cr
&& \prod_{v \not \in V} H^{i}(K_v,\widehat C) @<<<
\prod_{v \not \in U} H^{i}(K_v,\widehat C). 
\end{CD} 
$$
The right column is also exact by definition of $D^{i+1}(U,\widehat \C)$.
By diagram chasing, this yields an exact sequence 
\begin{equation} \label{vsequ}
\prod_{v \not \in U} H^i(K_v,\widehat C) \times \prod_{v \in U-V} 
H^i_{\nr}(K_v,\widehat C) \stackrel{s_V}{\to} H^{i+1}_c(V,\widehat \C) 
\to D^{i+1}(U,\widehat C).
\end{equation}
As seen before, the kernel of $s_V$ is the image of $H^i(U,\widehat \C)$, which 
implies that for $W \subset V$, the transition map 
$\ker s_W \to \ker s_V$ is surjective. The map 
$$\prod_{v \not \in U} H^i(K_v,\widehat C) \times \prod_{v \in U-W} 
H^i_{\nr}(K_v,\widehat C) \to \prod_{v \not \in U} H^i(K_v,\widehat C)
\times \prod_{v \in U-V} H^i_{\nr}(K_v,\widehat C)$$
is also obviously surjective. Thus taking projective limit 
over $V$ in (\ref{vsequ}) gives an exact sequence
$$\prod_{v \not \in U} H^i(K_v,\widehat C) \times \prod_{v \in U}
H^i_{\nr}(K_v,\widehat C) \to H^{1-i}(K,C)^* \to 
D^{i+1}(U,\widehat C)$$
(indeed recall that $H^{i+1}_c(V,\widehat \C) \cong H^{1-i}(V,C)^*$).
It remains to observe that for $i=2$, we have $D^3(U,\widehat \C) \subset
H^3(U,\widehat \C)=0$ (Proposition~\ref{finivanish2}).

\enddem

\begin{theo}[Poitou-Tate I] \label{spt1}
In the sequence 
\begin{equation} \label{ptdiag1}
\begin{CD}
0 @>>> H^{-1}(K,C) @>>> \prod^{'}_{v \in X^{(1)}} H^{-1}(K_v,C) @>>> 
H^2(K,\widehat C)^* \cr 
&&&&&& @VVV \cr
&& H^1(K,\widehat C)^* @<<< \prod^{'} _{v \in X^{(1)}} H^{0}(K_v,C)
@<<< H^0(K,C) \cr
&& @VVV \cr
&& H^1(K,C) @>>> \bigoplus_{v \in X^{(1)}} H^1(K_v,C) @>>> H^0(K,\widehat C)^* \cr
&&&&&& @VVV \cr
 0 @<<< H^{-1}(K,\widehat C)^* @<<< \bigoplus_{v \in X^{(1)}} H^2(K_v,C)
@<<< H^2(K,C) \cr
\end{CD}
\end{equation}

every sequence of three consecutive terms is exact, except the two ones 
respectively finishing with $H^0(K,C)$ and $H^1(K,C)$, which must be 
replaced with the following ``completed'' exact sequences: 
$$ \prod^{'}_{v \in X^{(1)}} H^{-1}(K_v,C)_{\wedge} \to H^2(K,\widehat C)^* 
\to H^0(K,C).$$
$$ \prod^{'}_{v \in X^{(1)}} H^{0}(K_v,C)_{\wedge} \to H^1(K,\widehat C)^* 
\to H^1(K,C).$$
\end{theo}

\dem{} Let $U \subset U_0$ be a non empty Zariski open subset with 
$U \neq X$. 

\smallskip

First take $i \in \{1, 2 \}$. By (\ref{suppcomp}), there is an exact sequence 
$$ H^i(U,\C) \to \bigoplus_{v \not \in U} H^i(K_v,C) \to 
H^{i+1}_c(U,\C)$$
and for $i=2$ the last map is surjective because $H^3(U,\C)=0$ 
by Proposition~\ref{finivanish1} a).
By Theorem~\ref{refin2} and  Proposition~\ref{injectcompl},
we have $$H^{i+1}_c(U,\C) \cong 
H^{1-i}(U,\widehat \C)^* \cong H^{1-i}(K,\widehat C)^* , $$
whence (by Lemma~\ref{vanishlocal}) for every non empty Zariski 
open subset $V \subset U$, a commutative diagram with exact lines 
$$
\begin{CD}
H^i(U,\C) @>>> \bigoplus_{v \not \in U} H^i(K_v,C) @>>> 
H^{1-i}(K,\widehat C)^* \cr 
@VVV @VVjV @VV=V \cr 
H^i(V,\C) @>>> \bigoplus_{v \not \in V} H^i(K_v,C) @>>> 
H^{1-i}(K,\widehat C)^*,
\end{CD}
$$ 
where $j$ is obtained by putting $0$ at the missing places (and the right
horizontal maps are surjective for $i=2$).
Therefore taking direct limit over $U$ in the first line of this diagram 
gives that the last two lines of (\ref{ptdiag1})
are exact.
The exactness of the first two lines of (\ref{ptdiag1}) comes from
Lemma~\ref{ptstart} after taking again direct limit over $U$. 

\smallskip

It remains to prove the exactness of the following three sequences:
$$ \prod^{'}_{v \in X^{(1)}} H^{-1}(K_v,C)_{\wedge} \to H^2(K,\widehat C)^* 
\to \Sha^0(C) \to 0.$$
$$ \prod^{'}_{v \in X^{(1)}} H^{0}(K_v,C)_{\wedge} \to H^1(K,\widehat C)^* 
\to \Sha^1(C) \to 0.$$
$$ \bigoplus_{v \in X^{(1)}} H^{1}(K_v,C) \to H^0(K,\widehat C)^* 
\to \Sha^2(C) \to 0.$$

We observe that for $0 \leq i \leq 2$, the following  
sequence is exact: 
$$0 \to \Sha^i(\widehat C) \to H^i(K,\widehat C) \stackrel{p_i}{\to}
\prod^{'}_{v \in X^{(1)}} H^i(K_v,\widehat C).$$
Set $A_i:={\rm Im} \, p_i$; we get exact sequences
(the maps being strict by Lemma~\ref{discretlem}) of 
Hausdorff, totally disconnected groups
$$0 \to \Sha^i(\widehat C) \to H^i(K,\widehat C) \to A_i \to 0.$$
$$0 \to A_i \to \prod^{'}_{v \in X^{(1)}} H^i(K_v,\widehat C),$$
where $A_i$ is equipped with the discrete topology. By \cite{hss}, Lemma~2.4
(where the groups are assumed to be locally compact, but the proof shows
that it sufficient to assume
that they have a basis of neighborhoods of zero consisting 
of open subgroups; this is the case for all groups considered here),
the duals of these exact sequences are also exact. 
Recall that for $i \in \{-1 , 0 \}$, the group
$H^i(K_v,C)_{\wedge}=\varprojlim_{n >0} (H^i(K_v,C)/n)$ is also the
completion $H^i(K_v,C)^{\wedge}$ of $H^i(K_v,C)$ with respect to the
open subgroups of finite index.
By \cite{demphd}, Th~3.1. and 3.3., the dual of 
$\prod^{'}_{v \in X^{(1)}} H^i(K_v,\widehat C)$ is 
$\prod^{'}_{v \in X^{(1)}} H^{1-i}(K_v,C)_{\wedge}$
for $1 \leq i \leq 2$, and
the dual of the group 
$\prod^{'}_{v \in X^{(1)}} H^0(K_v,\widehat C)=\prod_{v \in X^{(1)}} 
H^0(K_v,\widehat C)$ (cf. Remark~\ref{localisom}) is $\bigoplus_{v \in X^{(1)}} 
H^1(K_v,C)$. By Theorems~\ref{sha0theo} and \ref{sha12theo}, 
the dual of the finite group $\Sha^i(\widehat C)$ is 
$\Sha^{2-i}(C)$. This proves the result.

\enddem

\begin{theo}[Poitou-Tate II] \label{spt2}
In the sequence
\begin{equation} \label{ptdiag2}
\begin{CD}
0 @>>> H^{-1}(K,\widehat C) @>>> \prod_{v \in X^{(1)}}
H^{-1}(K_v,\widehat C) @>>> H^2(K,C)^* \cr 
&&&&&& @VVV \cr
&& H^1(K,C)^* @<<< \prod_{v \in X^{(1)}} H^{0}(K_v,\widehat C)
@<<< H^0(K,\widehat C) \cr
&& @VVV \cr
&& H^1(K,\widehat C) @>>> \prod^{'}_{v \in X^{(1)}} H^1(K_v,\widehat C) @>>>
H^0(K,C)^* \cr
&&&&&& @VVV \cr
 0 @<<< H^{-1}(K,C)^* @<<< \prod^{'}_{v \in X^{(1)}} H^2(K_v,\widehat C)
@<<< H^2(K,\widehat C) \cr
\end{CD}
\end{equation}

every sequence of three consecutive terms is exact, except the two ones
respectively finishing with $H^0(K,\widehat C)$ and $H^1(K,\widehat C)$,
which must be
replaced with the following completed exact sequences:
$$ \prod_{v \in X^{(1)}} H^{-1}(K_v,\widehat C)_{\wedge} \to H^2(K,C)^* 
\to H^0(K,\widehat C).$$
$$ \prod_{v \in X^{(1)}} H^{0}(K_v,\widehat C)_{\wedge}
\to H^1(K,C)^* 
\to H^1(K,\widehat C).$$
\end{theo}

\dem{} This is very similar to the proof of Theorem~\ref{spt1}. 
Let $U \subset U_0$ be a non empty open subset with $U \neq X$. 
For $i \in \{-1 , 0 \}$, (\ref{suppcomp}) and Proposition~\ref{finityp} a) 
give an exact sequence 
$$H^i(K,\widehat C) \to \prod_{v \not \in U} H^i(K_v,\widehat C) \to 
H^{i+1}_c(U,\widehat C),$$
such that the kernel of the first map is zero for $i=-1$, and this kernel 
is finite 
(by Lemma~\ref{d2lem}) for $i=0$. Applying Theorem~\ref{refin2} and 
taking projective limit over $U$, we obtain that the first two lines of 
(\ref{ptdiag2}) are exact.

\smallskip

Taking direct limit over $U$ in the exact sequences of Lemma~\ref{ptstart}
b) yields the exactness of the last line and of the sequence 
$$H^1(K,\widehat C) \to \prod'_{v \in X^{(1)}} H^1(K_v,\widehat C) \to 
H^0(K,C)^* \to \Sha^2(\widehat C).$$
because $\Sha^2(\widehat C) \cong \varinjlim_U D^2(U,\widehat \C)$ 
(cf. Theorem~\ref{sha0theo}).

\smallskip

Finally, dualizing the exact sequence of discrete groups 
$$0 \to \Sha^i(C) \to H^i(K,C) \to
\bigoplus_{v \in X^{(1)}} H^i(K_v,C)$$
for $i \in \{1 , 2 \}$  
gives the missing pieces of Theorem~\ref{spt2}, thanks to 
Theorem~\ref{sha12theo} and \cite{demphd}, Th. 3.1. and 3.3.

\enddem

\begin{rema} \label{variantors}
{\rm As the groups $H^1(K,\widehat C)$ and $H^2(K,\widehat C)$ are torsion, 
it is also possible to replace the last two lines of (\ref{ptdiag2}) by 
the following exact sequences
\begin{equation} \label{vari1}
0 \to \Sha^1(\widehat C) \to H^1(K,\widehat C) \to (\prod'_{v \in X^{(1)}}
H^1(K_v,\widehat C))_{\rm tors} \to (H^0(K,C)_{\wedge}) ^* \to
\Sha^2(\widehat C) \to 0.
\end{equation}

\begin{equation} \label{vari2}
0 \to \Sha^2(\widehat C) \to H^2(K,\widehat C) \to (\prod'_{v \in X^{(1)}} 
H^2(K_v,\widehat C))_{\rm tors} \to (H^{-1}(K,C)_{\wedge}) ^*  \to 0 .
\end{equation}
Indeed for $i \in \{ -1,0 \}$, the dual of $H^i(K,C)_{\wedge}=
\varprojlim_{n >0} (H^i(K,C)/n)$ (equipped with the inverse limit topology) 
is $$\varinjlim_{n >0} \, (H^i(K,C)/n)^*=\varinjlim_{n >0} \,  _n H^i(K,C)^*=
(H^i(K,C)^*)_{\rm tors},$$
which gives that (\ref{vari1}) and (\ref{vari2}) are exact, 
except that we don't have the surjectivity of 
$(H^0(K,C)_{\wedge}) ^* \to \Sha^2(\widehat C)$ yet. To see the latter, 
we first observe that for a $K$-torus $T$, the subgroup of divisible 
elements in $H^0(K,T)$ is trivial (indeed we may assume that $T$ is split 
and this is so for $K^*$ because $K$ is a global field); then the same 
property holds by d\'evissage (using (\ref{devistm})) for $H^0(K,C)$ because 
for a group of multiplicative type $M$, the group $H^1(K,M)$ is of finite 
exponent via Hilbert 90. Therefore the canonical map $H^0(K,C) \to 
H^0(K,C)_{\wedge}$ is injective, whence
an injection $\Sha^0(C) \hookrightarrow H^0(K,C)_{\wedge}$, whose dual 
$(H^0(K,C)_{\wedge})^* \to \Sha^2(\widehat C)$ (cf. Theorem~\ref{sha0theo})
is surjective, the group $\Sha^0(C)$ being finite and $H^0(K,C)_{\wedge}$
having a basis of neighborhoods of zero 
consisting of open subgroups (cf. \cite{hss}, Lemma~2.4).
}

\end{rema}

We can now prove the following variant of Theorem~\ref{spt1}:

\begin{theo}[Poitou-Tate, I'] \label{spt1bis}
There is an exact sequence 
\begin{equation} \label{ptdiag1bis}
\begin{CD}
0 @>>> H^{-1}(K,C)_{\wedge} @>>> [\prod^{'}_{v \in X^{(1)}}
H^{-1}(K_v,C)]_{\wedge} @>>>
H^2(K,\widehat C)^* \cr
&&&&&& @VVV \cr
&& H^1(K,\widehat C)^* @<<< [\prod^{'} _{v \in X^{(1)}} H^{0}(K_v,C)]_{\wedge}
@<<< H^0(K,C)_{\wedge} \cr
&& @VVV \cr
&& H^1(K,C) @>>> \bigoplus_{v \in X^{(1)}} H^1(K_v,C) @>>> H^0(K,\widehat C)^* \cr
&&&&&& @VVV \cr
 0 @<<< H^{-1}(K,\widehat C)^* @<<< \bigoplus_{v \in X^{(1)}} H^2(K_v,C)
@<<< H^2(K,C) \cr
\end{CD}
\end{equation}

\end{theo}

\dem{} Dualizing (\ref{vari1}) and (\ref{vari2}) 
yields the two exact sequences $$0 \to \Sha^0(C) \to 
H^0(K,C)_{\wedge} \to [\prod'_{v \in X^{(1)}} H^0(K_v,C)]_{\wedge} \to 
H^1(K,\widehat C)^* \to \Sha^1(C) \to 0.$$
$$0 \to H^{-1}(K,C)_{\wedge} \to [\prod'_{v \in X^{(1)}}
H^{-1}(K_v,C)]_{\wedge} \to H^2(K,\widehat C)^* \to \Sha^0(C) \to 0.$$
The other parts of the sequence follow from Theorem~\ref{spt1}.

\enddem

\begin{rema}
{\rm
A subtle point here is that for $i \in \{-1,0 \}$, 
there is a canonical injective map (which is induced by the isomorphism 
$\prod_v H^i(K_v,C)_{\wedge} \cong [\prod_v H^i(K_v,C)]_{\wedge}$):
$$\prod'_{v \in X^{(1)}} H^i(K_v,C)_{\wedge} \hookrightarrow 
[\prod'_{v \in X^{(1)}} H^i(K_v,C)]_{\wedge},$$
but this map is not surjective in general. For instance if $i=0$ and $C=\G$, 
we have $H^0(K_v,C) \cong \calo_v^* \times \Z$;
as $\bigoplus_{v \in X^{(1)}} \widehat \Z$ is smaller than 
$[\bigoplus_{v \in X^{(1)}} \Z]_{\wedge}$, the aforementioned map is not 
surjective. Observe that the natural map 
$H^i(K,C) \to \prod'_{v \in X^{(1)}} H^i(K_v,C)_{\wedge}$ {\it does not} 
in general extend to $H^i(K,C)_{\wedge}$, it is only defined on the 
bidual $H^i(K,C)^{**}$, which is smaller than $H^i(K,C)_{\wedge}$.
}
\end{rema}

\begin{rema} \label{II'}
{\rm Dualizing the exact sequence
$$0 \to \Sha^i(C) \to H^i(K,C) \to
\bigoplus_{v \in X^{(1)}} H^i(K_v,C) \to H^{1-i}(K,\widehat C)^*
\to \Sha^{i+1}(K,C) \to 0 $$
for $i \in \{1 , 2 \}$ also yields an exact Poitou-Tate sequence II', which 
is the same as (\ref{ptdiag2}) except that for $r \in \{-1, 0 \}$, the 
group $H^r(K,\widehat C)$ (resp. 
$\prod_{v \in X^{(1)}} H^r(K_v,\widehat C)$)
has to be replaced 
by $H^r(K,\widehat C)_{\wedge}$ (resp. 
by $(\prod_{v \in X^{(1)}} H^r(K_v,\widehat C))_{\wedge}=
\prod_{v \in X^{(1)}} H^r(K_v,\widehat C)_{\wedge}$).

}
\end{rema}

\begin{rema} \label{cdn3}
{\rm If we replace the function field $K$ by a number field, some 
results of this section still hold and some of them have to be 
modified. Namely: 

\begin{itemize}

\item Theorem~\ref{sha12theo} is unchanged (with the same proof), see
also Theorems~5.7 and 5.12 in \cite{demphd} (in the latter the 
assumption $\ker \rho$ finite is unnecessary). 

\item Lemma~\ref{injecth0} still holds. Therefore Theorem~\ref{sha0theo}
is also true: indeed since the pairing (\ref{d0pairing}) has divisible 
right kernel and trivial left kernel, taking the direct 
limit over $U$ (and using 
the facts that the sequence of finite groups $D^0(U,\C)$ stabilizes 
for $U$ sufficiently small) yields a pairing $\Sha^0(C) \times 
\Sha^2(\widehat C)$ with divisible right kernel and 
trivial left kernel. But it is known that $\Sha^2(\widehat C)$ is finite
(see \cite{demphd}, proof of Th. 5.14), whence the result (which extends 
\cite{demphd}, Th.~5.23).
Actually the image of $H^2_c(U,\widehat \C)$ into $H^2(K,\widehat C)$ is 
finite by d\'evissage thanks to exact triangle (\ref{devistmbis}):
indeed $H^3(K,\widehat T) \cong \bigoplus_{v \in \Omega_{\RR}}
H^3(K_v,\widehat T)$  is finite for a torus $T$, and 
$H^2_c(U,\widehat {\mathcal M})$
is also finite for a group of multiplicative type 
$\mathcal M$ because we already saw (cf. Remark~\ref{cdn1}
and Lemma~\ref{h1lem}, c) that 
this holds when $\mathcal M$ is a torus or a finite group.

\item In Lemma~\ref{discretlem}, one has to restrict to $i \geq 1$ for 
the assertion about $C$. The result about $\widehat C$ holds for 
an arbitrary $i$ (although the group $D^2(U,\widehat \C)$ might be infinite, 
we just saw that 
its image in $H^2(K,\widehat C)$ is finite, which is sufficient).

\item Lemma~\ref{ptstart} a) is not valid anymore (one has to complete
the first two terms in the exact sequences); the second exact sequence 
of b) still holds (same proof),
as does the first one except for the surjectivity of the last map,
which must be replaced by the exact sequence 
$$\prod_{v \not \in U} H^2(K_v,\widehat C) \times \prod_{v \in U} 
H^2_{\nr}(K_v,\widehat C) \to H^{-1}(K,C)^* \to D^3(U,\widehat \C) \to 0$$
because we lack the vanishing of $D^3(U,\widehat \C)$.
Also, since $D^2(U,\widehat \C)$ is in general infinite,
the proof of the exactness of 
$$H^2(U,\widehat \C) \to \prod_{v \not \in U} H^2(K_v,\widehat C)
\times \prod_{v \in U}
H^2_{\nr}(K_v,\widehat C) \to H^{-1}(K,C)^*$$
is a little bit more complicated (using exact triangle (\ref{devistmbis})
one reduces
to the case when $\C$ is quasi-isomorphic to $\mathcal M[1]$, where 
$\mathcal M$ is a group of multiplicative type; then one proceeds as in 
Lemma~\ref{ptstart} b), the group $D^2(U,\widehat {\mathcal M})$ being finite
because $H^2_c(U,\widehat {\mathcal M})$ is finite).

\item By the previous observations, the end of sequence (\ref{ptdiag1}) 
starting 
with $H^1(K,\widehat C)^* \to H^1(K,C) \to ...$ remains exact.
Theorem~\ref{spt2} is valid with one single slight complication in the 
proof: we do not know in general that $D^3(U,\widehat \C)=0$, 
but the direct limit over $U$ of the $D^3(U,\widehat \C)$ is 
$\Sha^3(\widehat C)$, which is zero. Theorem~\ref{spt1bis} is therefore also
unchanged, which extends \cite{demphd}, Theorems~6.1 and 6.3.

\end{itemize}

}
\end{rema}

\begin{rema}
{\rm In the case of one single torus $T$ with module of characters 
$\widehat T$, some of our results of section~\ref{two} and \ref{three} 
can be deduced from similar theorems on $1$-motives proven by 
Gonz\'alez-Avil\'es (\cite{gonz1}, Th. 6.6) and Gonz\'alez-Avil\'es/Tan 
(\cite{gonz2}, Th. 3.11).
}
\end{rema}

\end{document}